\documentclass[11pt]{article}
\usepackage{a4wide}
\usepackage{url}
\usepackage{booktabs}
\usepackage{longtable}
\usepackage{lscape}
\usepackage[table]{xcolor}
\usepackage{caption}
\usepackage{multirow}
\usepackage{amsmath, amssymb}
\usepackage{algorithm}
\usepackage{algorithmicx,algorithm, algpseudocode}
\usepackage{pstricks}
\usepackage{pst-all}
\allowdisplaybreaks
\usepackage[affil-it]{authblk}
\usepackage{tikz}
\usetikzlibrary{patterns}

\usepackage{subcaption}

\def\rev#1{{ #1}}

\newcommand{\itemset}{\mathcal{I}}

\newcommand\RotText[1]{\fontsize{9}{9}\selectfont \rotatebox[origin=c]{90}{\parbox{1cm}{\centering#1}}}

\begin{document}
\pagestyle{plain}
\date{}

\title{\rev{\Large \bf Exact Solution Techniques for\\ Two-dimensional Cutting and Packing}}
{\author{Manuel Iori$^{(1)}$, Vin\'icius L. de Lima$^{(2)}$, Silvano Martello$^{(3)}$,\\ Fl\'avio K. Miyazawa$^{(2)}$, Michele Monaci$^{(3)}$}
\affil{
(1) DISMI, University of Modena and Reggio Emilia (Italy)\\
(2) Institute of Computing, {University of Campinas (Brazil)}\\
(3) DEI ''Guglielmo Marconi'', University of Bologna (Italy)}
}
\date{}
\maketitle
\vspace*{-5ex}
\noindent

\begin{abstract}
We survey the main {formulations and} solution methods for two-dimensional orthogonal cutting and packing problems, where both items and bins are rectangles. We focus on exact methods and relaxations for the four main problems from the literature: {finding a packing with minimum height, packing the items into the minimum number of bins, finding a packing of maximum value, and determining the existence of a feasible packing.}
\end{abstract}


\section{Introduction}\label{sec:intro}

The number of publications on cutting and packing problems has been increasing considerably {in recent years.} In cutting problems, we are given a set of standardized stock units to be cut into smaller {\em{items}} so as to fulfill a given demand, while in packing problems a set of items has to be packed into one or more containers. These two classes of problems are strongly correlated, as packing an item into a container may be equivalent to cutting an item from a stock unit, and hence the same solution methods are often {adopted}. In the following, we will denote as {\em bins} both {the stock units and the containers}. In certain applications the unique container is a {\em strip} of (theoretically) infinite height.

Cutting and packing problems have been widely studied in the literature both for their theoretical interest and their many practical applications, in which they appear in a number of different variants.
Some problems consider one-dimensional items and bins,  whereas other problems refer to higher dimensions, like, e.g., the ({three-dimensional}) container loading problem (see Bortfeldt and W\"ascher \cite{BW13}). Items and bins may have rectangular or general (convex or non-convex) shape. Several practical constraints may also be part of the problems, as cargo stability in container loading problems or the need of producing guillotine patterns in cutting problems. In certain applications, cutting and packing problems are combined with other problems, for example in the integrated lot-sizing and cutting stock problems studied by Melega, Araujo, and Jans \cite{MAJ18} or in the routing problems with loading constraints considered by Iori and Martello \cite{IM10}. Packing problems also appear in many other fields, such as telecommunications \rev{(Lodi et al. \cite{LMMCLMEM11}, Martello \cite{M14})},
newspapers paging (Strecker and Hennig \cite{SH09}), production (Nesello et al. \cite{NDIS18}), {scheduling (Kwon and Lee \cite{KL15}), and maritime logistics (Xu and Lee \cite{XL18})}.
Most cutting and packing problems are ${\cal NP}$-hard and very challenging in practice. For this reason, sophisticated solution methods are needed for their solution, motivating the frequent  literature updating of this area of research.

{\subsubsection*{Books and surveys}}

\begin{figure}[h!]
\centering
\includegraphics[width=14cm]{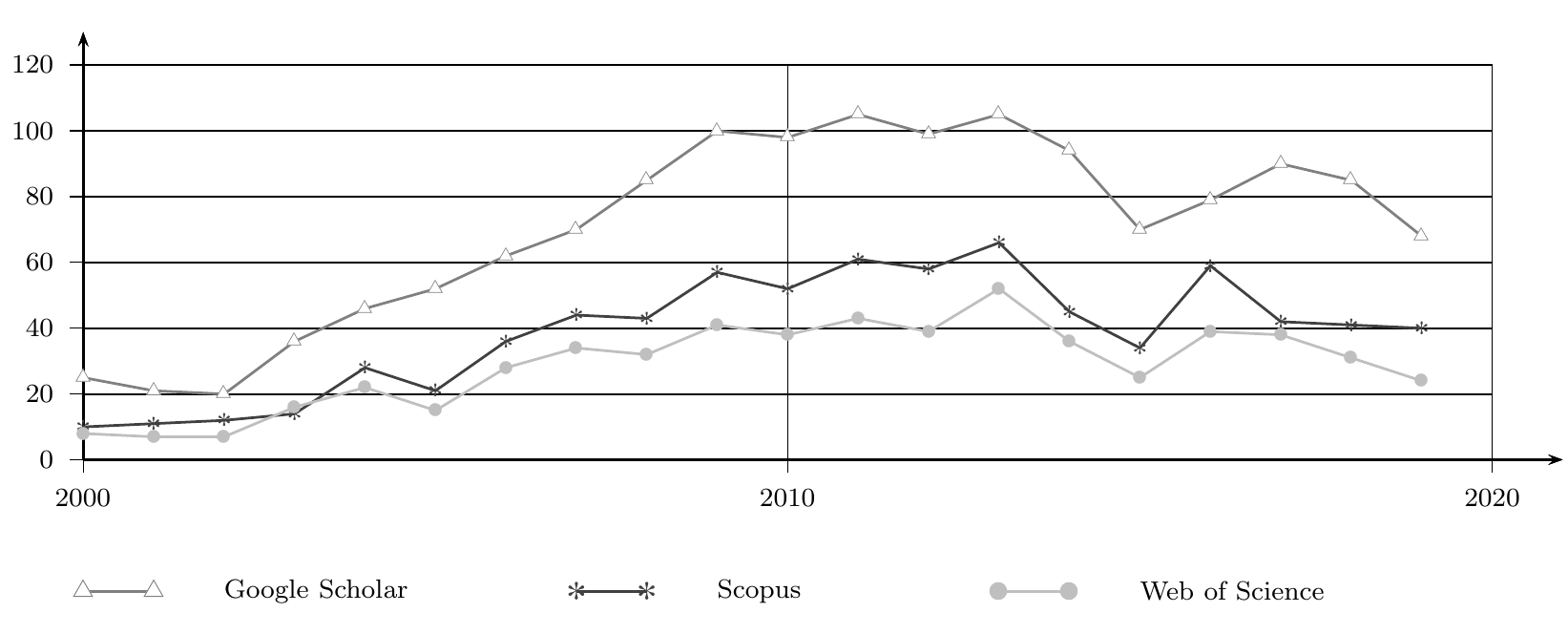}
\caption{\rev{Trend of number of publications on two-dimensional cutting and packing.}}
\label{figone}
\end{figure}

Several surveys have been dedicated to solution methods for cutting and packing problems. Specifically, referring to surveys published in the last two decades:
 \begin{itemize}
 \item {\em One-dimensional \rev{cutting and packing}}. Val{\'e}rio de Carvalho \cite{C02} reviews linear programming models for one-dimensional bin packing and cutting stock problems. \rev{Approximation algorithms for packing problems generally belong to two main categories: (i) {\em on-line} algorithms sequentially pack the items in the order encountered on input, without knowledge of items not yet packed; (ii) {\em off-line} algorithms
     have complete information about the item set, and may perform preprocessing, reordering, grouping, etc. before packing.
     Approximation results (for both on-line and off-line algorithms) have been reviewed by Coffman et al. \cite{CCGMV13} (Sections 3 and 4, respectively). A recent survey on mathematical models and exact algorithms for one-dimensional packing problems was presented} by Delorme, Iori, and Martello \cite{DIM16}, who also set up a library, the BPPLIB \cite{DIM18}, of computer codes and benchmark instances (see \url{http://or.dei.unibo.it/library/bpplib}). {The books by Martello and Toth \cite{MT90} and Kellerer, Pferschy, and Pisinger \cite{KPP04}
     present a comprehensive treatment on the knapsack problem and its variants but do not consider the corresponding two-dimensional versions};
 \item {\em Two-dimensional \rev{rectangular shape cutting and packing}}. After the classical 2002 surveys by Lodi, Martello, and Vigo \cite{LMV02} and Lodi, Martello, and Monaci \cite{LMM02},  a partial updating was presented
     by the same authors in 2014 \cite{LMMV14}. The following results {appeared} in the last decade. A survey on guillotine packing was produced by Ntene and van Vuuren \cite{NV09}.
     Silva, Oliveira, and W{\"a}scher \cite{SOW16} reviewed exact and heuristic algorithms for a particular problem (pallet loading) in which all items are identical, including a thorough analysis of benchmark instances and methodologies adopted in the literature for numerical experiments. Oliveira et al. \cite{ONJSC16} presented a review of {heuristic} algorithms for the strip packing problem. Christensen et al. \cite{CKPT17} proposed a survey {on approximation and on-line algorithms, also including} a list of open problems in this area. Recently, \rev{Russo et al. \cite{RBSS20} presented} a {survey on} relaxations for two-dimensional cutting problems with guillotine constraints and a categorization of the resulting bounds, {while Bezerra et al. \cite{BLOS19} reviewed models for the two-dimensional level strip packing problem};
 \rev{
 \item {\em Two-dimensional irregular shape cutting and packing.} To the best of our knowledge, the first survey dedicated to the packing of irregular shapes into rectangular containers was presented in the Nineties by Dowsland and Dowsland \cite{DD95}. More recently, two tutorials by Bennell and Oliveira \cite{BO08,BO09} reviewed the main geometric methodologies and algorithmic approaches for the heuristic solution of these problems. The latest survey on this area, presented by Leao et al. \cite{LTOCA20} in 2020, provides an extensive review of mathematical models for packing irregular shaped objects both in rectangular and irregular containers;
}
 \item {\em Multi-dimensional \rev{cutting and packing}}. In 1990, Dyckhoff \cite{D90} proposed a typology of cutting and packing problems (in one, two and three dimensions). About twenty years later, a successful improved typology  was proposed by W\"ascher, Hau$\beta$ner, and Schumann \cite{WHS07}.  \rev{More recently, Bortfeldt and W\"ascher \cite{BW13} considered the ({three-dimensional}) container loading problem and reviewed modeling approaches, as well as exact and heuristic algorithms.} A  recent survey on multi-dimensional packing problems was presented by Crainic, Perboli, and Tadei \cite{CPT12};
\item {\em Integrated variants.} {Reviews of algorithms for integrated routing and packing problems (with two- and three-dimensional packing constraints) have been presented by Iori and Martello \cite{IM10,IM13} and Pollaris et al. \cite{PBCJL15}.} Melega, de Araujo, and Jans \cite{MAJ18} recently reviewed integrated lot-sizing and cutting stock problems.
\end{itemize}
The \rev{constant interest in cutting and packing problems is shown by:}
{
\begin{itemize}
\item \rev{the trend of number of publications in the last 20 years, according to the major databases, shown in Figure \ref{figone}. In particular, the total number of publications in the last 20 years has been 1410, 774, and 575 according to Scholar, Scopus, and WoS, respectively;}
\item a number of special issues devoted to this area by international journals like, e.g., {\em INFOR} (see Martello \cite{M94}), {\em European Journal of Operational Research} (see Oliveira and W{\"a}scher \cite{OW07}), {\em International Journal of Production Economics} (see Bennell, Oliveira, and W\"ascher \cite{BOW13});
\item the recent comprehensive book by Scheithauer \cite{S18} dedicated to cutting and packing optimization;
\item the visual application for two-dimensional packing problems made available by Costa et al. \cite{CDIMM17};
\item the working group on cutting and packing (ESICUP) of the {\em Association of European Operational Research Societies} (EURO), see \url{https://www.euro-online.org/web/ewg/25/}.
\end{itemize}
}

{\subsubsection*{Contents}}
In this paper, we propose an extensive review of cutting and packing problems of two-dimensional rectangular items. We concentrate in particular on {\em orthogonal} cutting and packing, i.e., on the case in which
the items must be cut/packed with their edges parallel to those of the bin. {Some authors include the term ``geometric'' in the problem names.} The main problems that we address are:
{
 \begin{itemize}
 \item the two-dimensional strip packing problem: find a packing of minimum height into a single bin with fixed width;
 \item the two-dimensional bin packing problem: determine the minimum number of bins needed to pack the items;
 \item the two-dimensional knapsack problem: find a packing of maximum value {into a single bin};
 \item the two-dimensional orthogonal packing problem: find a feasible item packing (if any) into a single bin.
 \end{itemize}
}
We also consider some generalizations of these problems, such as the cutting stock problem, as well as relevant variants, e.g., the case in which the items can be rotated by 90 degrees or guillotine cuts are required.
We review the main relaxation methods, along with the main heuristics and exact methods.

The literature often considers formulations for two-dimensional cutting and packing problems which are based on a set of {points where an item may be packed or cut. We will discuss} the main methods for generating such sets of points. Another contribution of our work is the {review of preprocessing methods for cutting and packing problems, i.e., methods that are considerably fast, if compared to the overall solution time, and may reduce the size of the instances, thus} simplifying their resolution.

\rev{
The remainder of this survey is organized as follows. Section \ref{sec:problems} formally introduces the basic problems we address and the main variants considered in the literature. Section \ref{sec:preprocessing} reviews the principal techniques that can be used, in a preprocessing phase, to simplify a problem instance. Section \ref{sec:relaxations} discusses relaxation methods that provide valid bounds and, in some cases, are the starting point to derive a heuristic solution. Section \ref{sec:heuristics} briefly mentions the main heuristic techniques that are frequently embedded into exact approaches (although the scope of this survey is not to extensively review the huge literature on this topic). Section \ref{sec:ILP} examines mathematical models that explicitly require a solver, by classifying them according to their size (polynomial, pseudo-polynomial, or exponential). Section \ref{sec:enumeration} reviews enumeration schemes that do not explicitly make use of a solver (branch-and-bound, graph-based approaches, and constraint programming).}
\rev{Section \ref{sec:open_problems} provides pointers to relevant unsolved or challenging instances. Finally, Section \ref{sec:conclusions} discusses the state-of-the-art for each problem and presents some conclusions and future research directions.
}

\section{Problems and Definitions}\label{sec:problems}

In this section, we present a formal definition of the main problems and variants {we consider.} Given a bin and a set of items, a {\em packing} consists of a placement of the items such that the items {lie} completely inside the bin and there is no overlapping between any pair of items.

We are mainly interested in two-dimensional packing problems, in which the bins and the items are rectangles. Only orthogonal packings will be considered, i.e., we require that all edges of all items
are parallel to {the edges of the bin.} Unless otherwise specified, we assume that the items have {{\em fixed orientation}: as most of the} literature considers this case, we will only address as a variant the case in which the items can be rotated {by 90 degrees}.

A rectangular {\em two-dimensional bin} $\mathcal{B}$ is defined by its {\em width} $W\in \mathbb{Z}_+$ and {\em height} $H \in \mathbb{Z}_+$. We denote by $\itemset$ a set of rectangular two-dimensional {\em items}.
Each item $i \in \itemset$ has {\em width} $w_i \in \mathbb{Z}_+$ and {\em height} $h_i \in \mathbb{Z}_+$, such that $0 < w_i \leq W$ and $0 < h_i \leq H$. \rev{Unless otherwise specified, copies of identical items are treated as distinct items.}
A {\em packing} of $\itemset$ into $\mathcal{B}$ can be represented by a function $\mathcal{F}: \itemset \rightarrow \mathbb{Z}_+^2$ that maps each item $i \in \itemset$ to a pair $\mathcal{F}(i) = (x_i, y_i)$ representing the relative coordinates of the bottom-left corner of the item
with respect to the bottom-left corner of the bin. $\mathcal{F}$ is a {\em feasible packing} if
\begin{align}
&x_i \in \{0,...,W-w_i\} \text{ and } y_i \in \{0,...,H-h_i\} & (i \in \itemset) \label{constraint:bin_limits}\\*[1ex]
&[x_i,x_i+w_i)\cap [x_j,x_j+w_j) = \emptyset \text{ or }  [y_i,y_i+h_i)\cap [y_j,y_j+h_j) = \emptyset &(i, j \in \itemset, i \not = j).\label{constraint:no_overlapping}
\end{align}

In other words, the bin is seen on a Cartesian plane with its edges parallel to the $x$ and $y$ axes and its bottom-left corner on the origin. A packing defines, for every item $i$, the coordinates $(x_i, y_i)$ where its bottom-left corner is placed. Constraints \eqref{constraint:bin_limits} impose that each item is entirely inside the bin, while constraints \eqref{constraint:no_overlapping} forbid overlapping between any pair of items.

\subsection{Problems}\label{subsec:problems}
\begin{figure}[h!]
\rev{
\begin{picture}(50,00)(0,170)
	\resizebox{0.5\textwidth}{!}{
	\begin{tikzpicture}
	\draw[fill=white,thick] (9.75, 3.25) rectangle node{1} (10.25, 5.25); 	
	\draw[fill=white,thick] (8.50, 0.00) rectangle node{2} ( 9.50, 3.00);	
	\draw[fill=white,thick] (2.25, 0.00) rectangle node{3} ( 3.75, 3.00);	
	\draw[fill=white,thick] (4.00, 0.00) rectangle node{4} ( 5.50, 2.00);	
	\draw[fill=white,thick] (2.75, 3.25) rectangle node{5} ( 5.75, 5.25);	
	\draw[fill=white,thick] (0.00, 4.25) rectangle node{6} ( 2.50, 5.25);	
	\draw[fill=white,thick] (5.75, 1.25) rectangle node{7} ( 8.25, 2.25);	
	\draw[fill=white,thick] (0.00, 0.00) rectangle node{8} ( 2.00, 4.00);	
	\draw[fill=white,thick] (6.00, 3.25) rectangle node{9} ( 9.50, 5.25);	
	\draw[fill=white,thick] (5.75, 0.00) rectangle node{10}( 8.25, 1.00);	
	\draw[fill=white,thick] (5.10,-0.40) rectangle node{(a)}(5.10,-0.40);
	\end{tikzpicture}
	}
\end{picture}
%
\begin{picture}(50,0)(-275,252)
	\resizebox{0.25\textwidth}{!}{
	\begin{tikzpicture}
	\fill[fill=lightgray] (0,0) rectangle (5,9);
	\draw[very thick] (0,0) -- (5,0);
	\draw[very thick] (0,0) -- (0,9);
	\draw[very thick] (5,0) -- (5,9);
	\draw[dash pattern=on 6.3pt off 5pt, line width = 1.5pt] (0,9) -- (5,9);
	\draw[fill=white,thick] (3.50, 0.00) rectangle node{1} (4.00, 2.00);
	\draw[fill=white,thick] (4.00, 0.00) rectangle node{2} (5.00, 3.00);
	\draw[fill=white,thick] (0.00, 2.00) rectangle node{3} (1.50, 5.00);
	\draw[fill=white,thick] (1.50, 3.00) rectangle node{4} (3.00, 5.00);
	\draw[fill=white,thick] (0.00, 5.00) rectangle node{5} (3.00, 7.00);
	\draw[fill=white,thick] (1.50, 2.00) rectangle node{6} (4.00, 3.00);
	\draw[fill=white,thick] (2.50, 7.00) rectangle node{7} (5.00, 8.00);
	\draw[fill=white,thick] (3.00, 3.00) rectangle node{8} (5.00, 7.00);
	\draw[fill=white,thick] (0.00, 0.00) rectangle node{9} (3.50, 2.00);
	\draw[fill=white,thick] (0.00, 7.00) rectangle node{10}(2.50, 8.00);
	\draw[fill=white,thick] (2.50,-0.40) rectangle node{(b)}(2.50,-0.40);
	\end{tikzpicture}
	}
\end{picture}
%
\begin{picture}(50,0) (107,445)
	\resizebox{0.25\textwidth}{!}{
	\begin{tikzpicture}
	\fill[fill=lightgray] (0,0) rectangle (5,5);
	\draw[very thick] (0, 0) -- (5, 0);
	\draw[very thick] (0, 0) -- (0, 5);
	\draw[very thick] (5, 0) -- (5, 5);
	\draw[very thick] (0, 5) -- (5, 5);
	\fill[fill=lightgray] (0,5.5) rectangle (5,10.5);
	\draw[very thick] (0, 5.5) -- (5, 5.5);
	\draw[very thick] (0, 5.5) -- (0, 10.5);
	\draw[very thick] (5, 5.5) -- (5, 10.5);
	\draw[very thick] (0, 10.5) -- (5, 10.5);
	\draw[fill=white,thick] (0.00, 0.00) rectangle node{6} (2.5, 1.0);		
	\draw[fill=white,thick] (2.50, 0.00) rectangle node{7} (5.0, 1.0);		
	\draw[fill=white,thick] (0.00, 1.00) rectangle node{8} (2.0, 5.0);		
	\draw[fill=white,thick] (2.00, 1.00) rectangle node{5} (5.0, 3.0);		
	\draw[fill=white,thick] (2.00, 3.00) rectangle node{4} (3.5, 5.0);		
	\draw[fill=white,thick] (3.50, 3.00) rectangle node{1} (4.0, 5.0);		
	\draw[fill=white,thick] (0.00, 5.50) rectangle node{9} (3.5, 7.5);		
	\draw[fill=white,thick] (3.50, 5.50) rectangle node{3} (5.0, 8.5);		
	\draw[fill=white,thick] (0.00, 7.50) rectangle node{2} (1.0,10.5);		
	\draw[fill=white,thick] (1.00, 7.50) rectangle node{10}(3.5, 8.5);		
	\draw[fill=white,thick] (2.50,-0.40) rectangle node{(c)}(2.50,-0.40);
	\end{tikzpicture}
	}
\end{picture}
%
\begin{picture}(50,0) (-167.5,445)
	\resizebox{0.25\textwidth}{!}{
	\begin{tikzpicture}
	\fill[fill=lightgray] (0,0) rectangle (5,5);
	\draw[very thick] (0, 0) -- (5, 0);
	\draw[very thick] (0, 0) -- (0, 5);
	\draw[very thick] (5, 0) -- (5, 5);
	\draw[very thick] (0, 5) -- (5, 5);
	\draw[fill=white,thick] (4.5, 0) rectangle node{1} (5, 2);	
	\draw[fill=white,thick] (1.5, 2) rectangle node{2} (2.5, 5);	
	\draw[fill=white,thick] (0, 0) rectangle node{3} (1.5, 3);	
	\draw[fill=white,thick] (0, 3) rectangle node{4} (1.5, 5);	
	\draw[fill=white,thick] (1.5, 0) rectangle node{5} (4.5, 2);	
	\draw[fill=white,thick] (2.5, 2) rectangle node{6} (5, 3);	
	\draw[fill=white,thick] (2.5, 3) rectangle node{7} (5, 4);	
	\draw[fill=white,thick] (2.5, 4) rectangle node{10} (5, 5);	
	\draw[fill=white,thick] (2.50,-0.40) rectangle node{(d)}(2.50,-0.40);
	\end{tikzpicture}
	}
\end{picture}
}

\vspace*{95ex}
\caption{
(a) {set of items};
(b) an optimal 2D-SPP solution;
(c) an optimal 2D-BPP solution;
(d) an optimal 2D-KP solution (in case items profits correspond to their areas).
}
\label{fig:example_problems}
\end{figure}
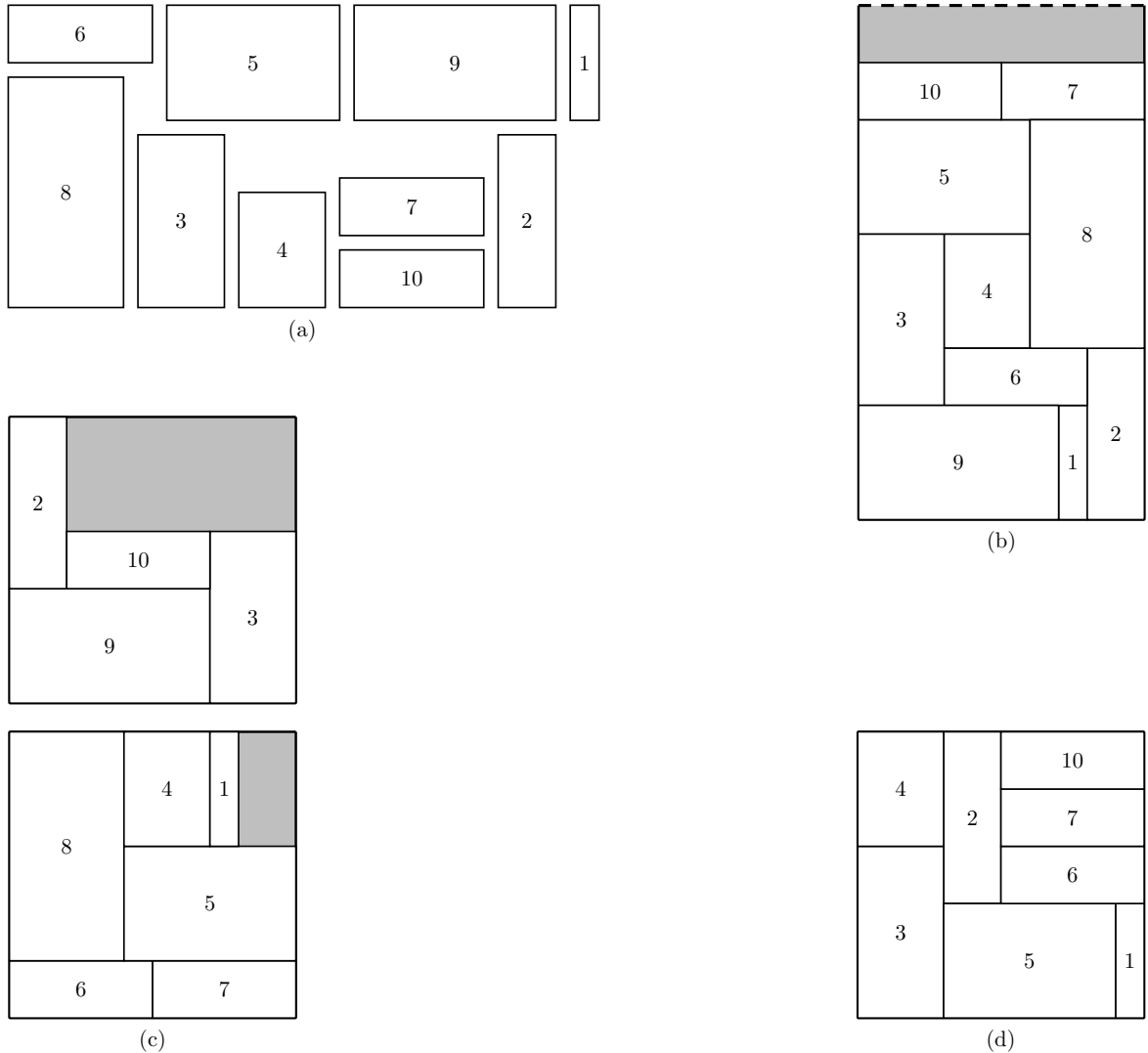

In this section, we define the main two-dimensional orthogonal packing problems {we consider.
As mentioned in Section \ref{sec:intro},} some typologies have been proposed in the literature to classify the great variety of possible packing problems. As we will mainly deal with the four problems defined below, we preferred to use the simple names we provide. \rev{For the sake of completeness, their definition is followed, in Section \ref{sec:typology}, by the corresponding notation in the main typologies.}  \rev{A simple visual example of a set \rev{$\itemset$ containing 10 items} is given in Figure \ref{fig:example_problems}(a), and then used to clarify the main problem variants below.}

\subsubsection*{Two-Dimensional Strip Packing Problem}
{We are given a single bin $\mathcal{B}$ having fixed width $W$ and infinite height, usually called a {\em strip}. The {\em Two-Dimensional Strip Packing Problem} (2D-SPP) asks for} a feasible packing of $\itemset$ into $\mathcal{B}$ such that the height at which the strip is used (i.e., the height of the topmost edge of an item) \rev{is minimized}.
\rev{Figure \ref{fig:example_problems}(b) shows a minimum height arrangement of the items of Figure \ref{fig:example_problems}(a) into a strip.}

\subsubsection*{Two-Dimensional Bin Packing Problem}
{In this case we have an unlimited number of identical finite bins $\mathcal{B}$ having width $W$ and height $H$. The} {\em Two-Dimensional Bin Packing Problem} (2D-BPP) requires to determine a partition of $\itemset$ \rev{into the minimum} number of subsets, such that each subset can be feasibly packed into a bin. A generalization of the problem is the {\em Two-Dimensional Cutting Stock Problem} (2D-CSP), in which one is asked to pack a given number $d_i$ ({\em demand}) of each item $i \in \itemset$.
Both problems generalize their one-dimensional version, the well-known ({one-dimensional})  {\em Bin Packing Problem} (1D-BPP) and {\em Cutting Stock Problem} (1D-CSP), in which the items are segments of size $w_i$
($i \in \itemset$) and the bins are segments of size $W$ {({\em capacity})}. An extensive literature exists on these problems, see \cite{DIM16} for a recent survey.
\rev{Figure \ref{fig:example_problems}(c) \rev{shows an optimal 2D-BPP solution, in which the items of Figure \ref{fig:example_problems}(a)} are packed into two separate bins.}

\subsubsection*{Two-Dimensional Knapsack Problem}
The problems listed so far ask for a feasible packing of {\em all} {items of $\itemset$}. Assume now that every item $i \in \itemset$ has an associated value ({\em profit}) $v_i \in \mathbb{Z}_+$. The {\em Two-Dimensional Knapsack Problem} (2D-KP) requires to determine a subset of items $\itemset' \subseteq \itemset$ such {that: (i) there exists a feasible packing of $\itemset'$ into a single bin $\mathcal{B}$; and (ii) the} corresponding total profit, $\sum_{i \in \itemset'} v_i$, \rev{is maximized}. A special case of the 2D-KP, the {\em Two-Dimensional Rectangular Packing Problem}, arises when the profit of each item is equal to {its area (i.e., \rev{$v_i = w_i h_i$} $\forall$ $i \in \itemset$)}.
\rev{An example of an optimal \rev{two-dimensional rectangular packing problem solution for the items of Figure \ref{fig:example_problems}(a)} is given in Figure \ref{fig:example_problems}(d).}

\subsubsection*{Two-Dimensional Orthogonal Packing Problem}
{While the above problems are in {\em optimization version}, the} {\em Two-Dimensional Orthogonal Packing Problem} (2D-OPP) is to decide if there exists a feasible packing of a given set $\itemset$ of items into a single bin $\mathcal{B}$. {Although this problem is in {\em decision version}, i.e., it just asks for a `yes/no' answer, in practical applications one is generally requested to also produce the specific packing, if any. The 2D-OPP  appears as a subproblem in the optimization version of a number of packing problems.}

{
\subsection{Typologies}\label{sec:typology}
}
As previously mentioned, a number of typologies have been introduced in the literature. \rev{As the reader could encounter the same problems considered in the present survey but identified in a different way, we next provide the notations adopted by the most common typologies from the literature}.

{According to the typology proposed by Dyckhoff \cite{D90} in 1990, the 2D-BPP, the 2D-CSP and the 2D-KP are denoted as 2/V/I/M, 2/V/I/R, and 2/B/O/, respectively. The 2D-SPP and the 2D-OPP have not been formally defined in this typology.}

{In 1999 Lodi, Martello, and Vigo \cite{LMV99a} used a three-field typology (later extended by Martello, Monaci, and Vigo \cite{MMV03}). The 2D-SPP, the 2D-BPP, and the 2D-KP are denoted as 2SP$|$O$|$F, 2BP$|$O$|$F, and 2KP$|$O$|$F, respectively, while the 2D-OPP is not classified. The second field of this notation also covers the variant (see below) in which orthogonal rotation of the items is allowed (`R' instead of `O'), while the third field can handle the variant in which guillotine cuts are required (`G' instead of `F').}

{More recently (2007), W\"ascher, Hau$\beta$ner, and Schumann \cite{WHS07} proposed a successful typology, partially based on Dyckhoff's original ideas. According to it,} the 2D-SPP is {denoted as the ``(Two-dimensional) Open Dimension Problem'' (ODP), the 2D-BPP as the ``Single Bin Size Bin Packing Problem'' (SBSBPP), and the 2D-CSP as the ``Single Stock Size Cutting Stock Problem'' (SSSCSP). The 2D-KP is denoted either as the ``Single Knapsack Problem'' (SKP), or as the ``Single Large Object Placement Problem'' (SLOPP) in the generalization in which
each item $i \in \itemset$ is available in $d_i$ copies. The 2D-OPP is not explicitly classified in \cite{WHS07}, but it can be interpreted as the recognition version of the
2D-BPP.}

{\subsection{Complexity}}

{
The {\em one-dimensional} version of the 2D-KP is the well-known {\em knapsack problem} (1D-KP) in which each item has a {\em weight} $w_i$ and the bin ({\em knapsack}) has capacity $W$. The 1D-KP can be solved in pseudo-polynomial time through dynamic programming (see, e.g., Martello and Toth \cite{MT90} or Kellerer, Pferschy, and Pisinger \cite{KPP04}). Instead, none of the above two-dimensional problems admits a pseudo-polynomial time algorithm, unless ${\cal P} = {\cal NP}$. Recall indeed that the 1D-BPP is strongly ${\cal NP}$-hard (see Garey and Johnson \cite{GJ79}). Given an instance of the 1D-BPP, define two-dimensional items having width $w_i$ and
height $h_i = 1$ ($i \in \itemset$). Then:
\begin{itemize}
\item the solution of a 2D-SPP instance with strip width $W$ solves the 1D-BPP instance;
\item the solution of a 2D-BPP instance with bins of width $W$ and height $1$ solves the 1D-BPP instance;
\item associate a profit $v_i = 1$ to each item $i \in \itemset$. {The minimum value $\bar{H}$ for which the optimal solution of a 2D-KP
instance with a bin of width $W$ and height $\bar{H}$ has value $n$ gives the optimal 1D-BPP solution value;}
\item the solution of a 2D-OPP instance with a bin of width $W$ and height $k$ {answers the decision version of the 1D-BPP: can the items} of $\itemset$ be packed into $k$ bins?
\end{itemize}
It follows that the existence of a polynomial-time algorithm for any of the four problems above would imply a polynomial-time algorithm for the 1D-BPP, i.e., the 2D-SPP, the 2D-BPP, and the 2D-KP are strongly ${\cal NP}$-hard, and the 2D-OPP is strongly ${\cal NP}$-complete.
}

\subsection{Variants}\label{sec:variants}
Motivated by practical applications, a number of variants of two-dimensional packing problems has been considered in the literature (see, e.g., Lodi, Martello, and Vigo \cite{LMV99a}, Pisinger and Sigurd \cite{PS05}, and W\"ascher, Hau$\beta$ner, and Schumann \cite{WHS07}).

\subsubsection*{Orthogonal Rotation}
While in the basic problems the items have fixed orientation, relevant variants allow item {\em rotation} by 90 degrees. {For example, while in the cutting of corrugated or decorated stock units rotation is forbidden, when the surfaces are uniform it may be feasible to rotate the items in order to produce better (more dense) packings. Note however that this also leads to more complex problems as it increases the number of decisions to be taken, and hence the corresponding models usually involve a higher number of variables and constraints.} \rev{Indeed, in most approaches, item rotation is typically handled either by adding, for each item,  a ``rotate-or-not'' binary decision (like, e.g., in Jakobs \cite{J96}), or by creating a companion (rotated) copy of each item and forbidding that both are selected (like, e.g., in Lodi and Monaci \cite{LM03}).}

\subsubsection*{Guillotine Cuts}
\begin{figure}[b]
\centering
\begin{subfigure}[htpb]{0.32\textwidth}\centering
\resizebox{0.8\textwidth}{!}{
	\begin{tikzpicture}
	\draw[fill=lightgray,ultra thick] (0,0) rectangle (5,5);

	\draw[fill=white,very thick] (0,0) rectangle node{} (3.5,2);
	\draw[fill=white,very thick] (3.5,0) rectangle node{} (5,2.5);
	\draw[fill=white,very thick] (0,2) rectangle node{} (2,5);
	\draw[fill=white,very thick] (2,2.5) rectangle node{} (5,5);

	\end{tikzpicture}
}
\caption{\label{fig:non-guillotine}}
\end{subfigure}
\begin{subfigure}[htpb]{0.32\textwidth}\centering
\resizebox{0.8\textwidth}{!}{
	\begin{tikzpicture}
	\draw[fill=lightgray,ultra thick] (0,0) rectangle (5,5);

	\draw[fill=white,very thick] (0,0) rectangle node{} (3.5,2);
	\draw[fill=white,very thick] (0,2) rectangle node{} (2,5);
	\draw[fill=white,very thick] (2,2) rectangle node{} (5,4.5);
	\draw[fill=white,very thick] (3.5,0) rectangle node{} (5,1);

	\end{tikzpicture}
}
\caption{\label{fig:2staged-guillotine}}
\end{subfigure}
\begin{subfigure}[htpb]{0.32\textwidth}\centering
\resizebox{0.8\textwidth}{!}{
	\begin{tikzpicture}
	\draw[fill=lightgray,ultra thick] (0,0) rectangle (5,5);

	\draw[fill=white,very thick] (0,0) rectangle node{} (3.5,2);
	\draw[fill=white,very thick] (0,2) rectangle node{} (2,5);
	\draw[fill=white,very thick] (2,2) rectangle node{} (5,4.5);
	\draw[fill=white,very thick] (3.5,0) rectangle node{} (5,1);
	\draw[fill=white,very thick] (3.5,1) rectangle node{} (4.5,2);
	\draw[fill=white,very thick] (2,4.5) rectangle node{} (4,5);

	\end{tikzpicture}
}
\caption{\label{fig:3staged-guillotine}}
\end{subfigure}
\caption{ (a) non-guillotine pattern; (b) 2-staged guillotine pattern with trimming; (c) 3-staged guillotine pattern with trimming. }
\label{fig:guillotine}
\end{figure}
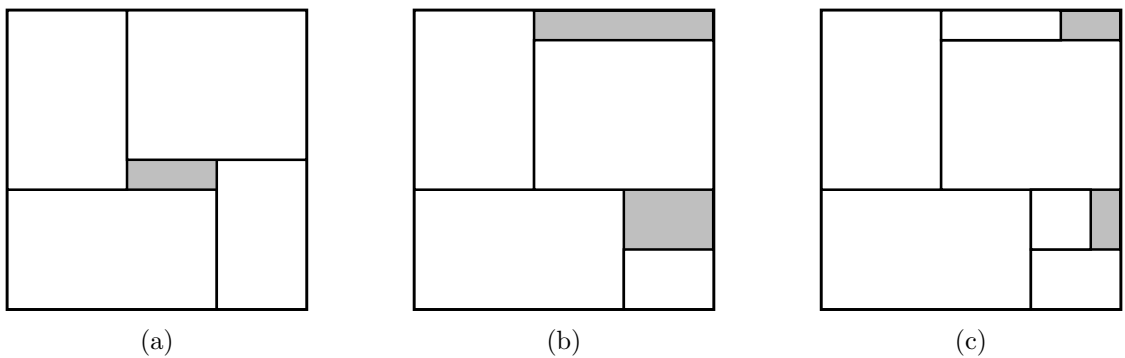

Normally, automatic cutting machines can only produce the items through a sequence of {\em guillotine cuts}, i.e., edge-to-edge cuts parallel to the edges of the bin. {Imposing such constraint may lead to worse solutions as not all item sets allow guillotine patterns in a bin (see Figure \ref{fig:guillotine}(a))}. Frequently, the machines are restricted to only
alternate horizontal and vertical cuts, possibly with a hard limit $k$ on the number of cuts per bin ({\em k-staged} problems). In most applications $k$ is two or three, possibly allowing an extra cut (called \emph{trimming}) to separate an item from waste {(see Figure \ref{fig:guillotine}(b) and (c))}. When $k = 2$, the problem is frequently called {\em level packing}, as it can be seen as the problem of packing the items side-by-side on horizontal shelves having width equal to that of the bin/strip and height coinciding with the tallest packed item {(see again Figure \ref{fig:guillotine}(b))}. \rev{When $k=3$, Puchinger and Raidl \cite{PR07} distinguish between the case where the latter condition holds ({\em restricted case}) and the one where the height of a shelf is not necessarily given by its highest item ({\em unrestricted case}).}
\rev{The presence of guillotine constraints consistently affects the combinatorial structure of the problem, and hence the solution techniques. Indeed, for such cases,
it is common to adopt techniques based on column generation and/or dynamic programming, that would not be suitable for standard non-guillotine problems}.

\subsubsection*{Variable-Sized Bins}
The bin packing and cutting stock problems consider an unlimited number of identical bins. The {\em variable-sized} generalization of these problems deals instead with different types of bins, each
having a specific size (width and height), cost, and availability. The problem is then to pack {all the items} at minimum cost.
\rev{The addition of variable sizes is usually handled in the solution techniques by considering, e.g., the different sizes of the bins when creating new nodes in enumerative approaches, or when pricing variables in branch-and-price methods (as, e.g., in Pisinger and Sigurd \cite{PS05}). This comes at the cost of an additional computational effort, and indeed, from a practical point of view, the solution of variable sized problems is
typically more challenging.}

\subsubsection*{Loading and Unloading Constraints}
Loading and/or unloading constraints often arise in applications where goods pertaining to different customers have to be loaded into the same {bin, which represents, e.g., the loading area of a truck}
\rev{(see, e.g., Gendreau et al. \cite{GILM06})}. In these cases, the vehicles have a loading/unloading orientation, and the sequence of visit to the customers must be such that each item may be moved in/out without moving any other item.
\rev{These constraints can be directly included in mathematical models where the variables indicate the position of the items (see Section \ref{sec:pseudo}) as the packing of an item restricts the range of positions that can be taken by other items (as, e.g., in C{\^o}t{\'e}, Gendreau and Potvin \cite{CGP14}). They can also be imposed within branch-and-bound algorithms, where additional fathoming criteria can be devised to reduce the number of decision nodes (as, e.g., in Iori, Salazar Gonz{\'a}lez and Vigo \cite{ISV07}).}\\

\rev{The first two variants were already discussed in the Sixties by Gilmore and Gomory \cite{GG65}. In the Eighties, Friesen and Langston \cite {FL1986} introduced (for the one-dimensional case) the variable-size variant. The studies on loading and unloading constraints started in the Noughties, see Gendreau et al. \cite{GILM06}.}

\section{Sets of Points and Preprocessing Techniques}\label{sec:preprocessing}

We start this section by reviewing techniques based on sets of points, that can be used in both {{\em Integer Linear Programming} (ILP) models and solution methods reviewed in the following sections.} We then examine some preprocessing techniques, used to decrease the size of {a given instance.}

\subsection{Sets of Points} \label{sec:points}

Several authors in the literature use {\em sets of points} to represent the possible positions where an item can be packed. According to the definitions given in Section \ref{sec:problems}, packing an item at a point
$p = (x,y)$ means to allocate the item into the bin with its bottom-left corner in the position identified by $p$.

As we assume that items and bins have integer sizes, only considering integer point coordinates does not affect optimality. Trivial sets of points for an item $i$ are thus given by coordinate sets $X_i = \{0,\dots,W-w_i\}$ for the $x$-axis and $Y_i = \{0,\dots,H-h_i\}$ for the $y$-axis.

It is clear that the smaller the number of points the smaller the search domain. We review in the following the main techniques that have been proposed in the literature for generating reduced sets of points by preserving optimality.

Already in the Seventies, Herz \cite{H72} and Christofides and Whitlock \cite{CW77} independently observed that any feasible packing pattern can be transformed into an equivalent one ({\em normal pattern}) in which the items are shifted to the bottom and left as much as possible. In a normal pattern, the left and bottom edges of each item touch another item or the bin. This leads to the definition of the set of ``normal'' $x$-coordinates
\begin{equation}
\label{formula:normal_patterns}
\mathcal{N}^{\rm x} = \{ x \,|\, x < W,\, x = \sum_{i \in \itemset} w_i \zeta_i,\, \zeta_i \in \{0,1\} (i \in \itemset) \},
\end{equation}
\noindent
{defining the only $x$-coordinates at which an item can be packed.} A similar definition holds for the set $\mathcal{N}^{\rm y}$ of normal $y$-coordinates.

Boschetti, Hadjiconstantinou, and Mingozzi \cite{BHM02} proposed to separately compute a set of normal patterns for each {item $k \in \itemset$ by excluding those patterns that include $k$ itself. The resulting
set $\mathcal{N}_k^{\rm x}$ ($k \in \itemset$) is obtained from \eqref{formula:normal_patterns} by replacing $\itemset$ with $\itemset \setminus \{k\}$. (Similarly for $\mathcal{N}_k^{\rm y}$.)}

{Terno, Lindemann, and Scheithauer \cite{TLS87} proposed the use of \emph{reduced raster points}, obtained by removing redundant positions from the set of normal patterns. The idea is that, given two coordinates $p, q \in \mathcal{N}^{\rm x}$ ($p<q$), if every possible combination of items that can be packed to the right of $p$ can also be packed to the right of $q$, then $p$ can be removed {from $\mathcal{N}^{\rm x}$.}

{Recently, C{\^o}t{\'e} and Iori \cite{CI18} proposed a new set of patterns {called {\em meet-in-the-middle}, generated} by first defining a threshold $t \in \{1, 2, \dots, W\}$ and then left-aligning patterns that are to the left of $t$ and right-aligning those that are to its right. The resulting set of patterns is proved to be never larger than $\mathcal{N}^{\rm x}$, and in practice is usually much smaller.}

\rev{The fact that normal patterns preserve optimality was proved in the seminal paper by Herz \cite{H72}. Proofs were later provided for other types of patterns like, e.g., in C{\^o}t{\'e} and Iori \cite{CI18}. All such proofs show that any solution in which the items are packed without restrictions in the continuous two-dimensional space can be transformed, through simple translations, into a better or equivalent solution in which the items are packed according to the  patterns. The same arguments apply to problem variants like, e.g., those involving guillotine cuts or several item copies or item rotation.}

\subsection{Preprocessing Techniques}

Two main approaches for preprocessing the instances by preserving optimality can be found in the literature: methods that fix some decisions and methods that modify the input parameters. These techniques are useful for improving bounds based on the area \rev{of the items or of the bins. We describe in the following one method of the former type and two of the latter. Although developed for specific problems, the methods of the latter type can be directly used for all problems. The extension of methods of the former type is instead not straightforward.}

\rev{All methods will be presented referring to widths. Those developed for problems on bins
can be identically applied using heights instead of widths, while the same does not hold for problems on strips. The sequence of presentation follows the order in which it is advisable to execute these methods in order to effectively reduce the given instance.}

\subsubsection*{Packing a Subset of Items}
This reduction was originally proposed by Martello, Monaci, and Vigo \cite{MMV03} for the 2D-SPP. Suppose that the items are sorted by non-increasing width. The method can be applied if $w_1 > W/2$. Let
$B = \{ \ i \in \itemset \ | \ w_i = w_1 \ \}$. Observe that these items cannot be packed side-by-side, so we can pile them aligned to the left edge of the strip. Now define the set of those items that can be packed side-by-side with an item of $B$, namely $S = \{ \ j \ | \ w_j \leq W - w_1 \ \}$. If there exists a feasible packing of all the items of $S$ into the empty right part of the strip, with overall height not greater than $\sum_{i \in B} h_i$,
then this packing can be optimally fixed on the bottom of the strip, and the process can be iterated on the reduced instance. When the required packing is not found, the method finds the first item $\ell \not \in B$.
{If $w_\ell > W/2$,} the current set $B$ is updated by adding the items of width $w_\ell$, the current set $S$ is updated accordingly, and a new attempt is performed.

This reduction was later used within other solution methods for the 2D-SPP (see, e.g., Alvarez-Valdes, Parre{\~n}o, and Tamarit \cite{APT09}, Boschetti and Montaletti \cite{BM10} and \rev{C{\^o}t{\'e}, Dell'Amico, and Iori \cite{CDI14})} and extended to other packing problems, among which the 2D-BPP (see, e.g., Carlier, Clautiaux, and Moukrim \cite{CCM07}).

\subsubsection*{Shrinking the Size of the Bin}
Alvarez-Valdes, Parre{\~n}o, and Tamarit \cite{APT09} proposed to solve a (one-dimensional) \rev{{\em subset-sum problem} (find a subset of a set of given integers whose sum is closest to, without exceeding, a prefixed threshold)}
to determine the maximum value $\overline W${$\leq W$} such that there exists a sum of item widths equal to $\overline W$.
If $\overline W < W$ then the width of the bin/strip can be reduced to $\overline W$.

\subsubsection*{Lifting the Size of the Items}
A similar approach had been used by Boschetti, Hadjiconstantinou, and Mingozzi \cite{BHM02} to increase item widths. They proposed to solve, for each item $i \in \itemset$, a subset-sum problem to determine the maximum value $\overline W_i \leq W - w_i$ such that there exists a set of items in $\itemset \setminus \{i\}$ with total width equal to $\overline W_i$.
If $w_i + {\overline W}_i < W$ then the width of $i$ can be increased by $W - {\overline W}_i$.
Carlier, Clautiaux, and Moukrim \cite{CCM07} defined similar methods for updating the size of the items by removing small items and increasing the width of large items.

\section{Relaxations}\label{sec:relaxations}

Several relaxation methods for two-dimensional packing problems have been proposed in the literature. They are used within exact algorithms and, in some cases, as a base to construct a heuristic solution.
{Obviously,} relaxations provide lower bounds for the 2D-SPP and the 2D-BPP and upper bounds for the 2D-KP. As the 2D-OPP is a decision problem, relaxations can be used to prove that the required packing does not exist.

\subsection{Continuous Relaxation}
Splitting each item $i$ into $w_i \times h_i$ unit squares produces the most immediate relaxation for all considered problems. For the 2D-OPP, the 2D-SPP and the 2D-BPP, the {\em continuous lower bound} is then
\begin{equation}
\label{formula:continuous_lb}
\left \lceil \sum_{i\in \itemset } \frac{w_i h_i}{WH} \right \rceil
\end{equation}
(computed with $H = 1$ for the 2D-SPP). For the 2D-KP, the {{\em continuous upper bound}} is obtained, following Dantzig \cite{D57} by: (i) sorting the items by non-increasing ratios $v_i/(w_i h_i)$; (ii) finding the first item $s$ such that $\sum_{i=1}^s w_i h_i > W H$ and defining $c = WH - \sum_{i=1}^{s-1} w_i h_i$; (iii) computing the upper bound as
\begin{equation}
\label{formula:continuous_ub}
\left \lfloor \sum_{i=1}^{s-1} v_i + \frac{c}{w_s h_s} v_s \right \rfloor.
\end{equation}
Efficient implementations allow all these bounds to be computed in linear time.
\rev{Because of its simplicity, the continuous relaxation  has been used in almost all works on two-dimensional packing problems.}

{\subsection{Combinatorial Bounds}}
Martello and Vigo \cite{MV98} extended to the 2D-BPP {some} lower bounds proposed by Martello and Toth \cite{MT90b} for the 1D-BPP and by
Dell'Amico and Martello \cite{DM95} for the $P||C_{\max}$, a parallel machine scheduling problem that is strictly related to the 1D-BPP. \rev{The idea is to identify, for a given parameter $p$:
(i) two item sets, say $J_1(p)$ and $J_2(p)$, such that no two items of $J_1(p) \cup J_2(p)$ may be packed into the same bin; (ii) a third set, $J_3(p)$, of items that cannot be packed into a bin used for an item of $J_1(p)$.
A valid lower bound is then $L(p) = |J_1(p)\cup J_2(p)| + L_{23}(p)$, where $L_{23}(p)$ denotes a lower bound on the number of additional bins needed for the items of $J_3(p)$. It is shown that
the overall lower bound, $\max_{p} \{L(p)\}$ can be computed in $O(n^2)$ time.}

Boschetti and Mingozzi \cite{BM03A,BM03B} improved these bounds and extended them to the variant in which orthogonal rotation of the items is allowed.
Similarly, relaxations for the 2D-KP can be obtained from induced 1D-KP instances where each item has a weight equal to its area and the knapsack has  capacity \rev{$W H$.}

\rev{Although combinatorial bounds have a low computational cost, they  can be quite effective, especially in branch-and-bound approaches (see, e.g., Martello, Monaci, and Vigo \cite{MMV03} and Alvarez-Valdes, Parre{\~n}o and Tamarit \cite{APT09}).}


\subsection{Linear Relaxation and Column Generation}\label{sec:column}
Already in 1965, the seminal paper by Gilmore and Gomory \cite{GG65} introduced a column generation algorithm for {two-dimensional} packing problems. The algorithm is based on the linear relaxation of a mathematical formulation (the so-called {\em set covering model}) that they had developed for the 1D-CSP (see Section \ref{subsec:problems}). The model has a variable for each possible combination of items that {fits into a} single bin ({\em pattern}). Since the resulting {\em Linear Program} (LP) has a huge (exponential) number of \rev{{\em columns} (variables)}, the {\em Column Generation} method starts with a restricted model that only includes a subset of columns and iteratively adds columns that may improve the current solution.
\rev{This relaxation has been used in branch-and-price algorithms for the 2D-BPP and related variants (see, e.g., Pisinger and Sigurd \cite{PS05,PS07}), and in problems involving guillotine cuts (see, e.g., Belov and Scheithauer \cite{BS06}, Bettinelli, Ceselli, and Righini \cite{BCR08}, and Cintra et al. \cite{CMWX08}).}

\subsection{Dual Feasible Functions}

{Dual feasible functions were originally introduced by Johnson \cite{J73} for a normalized version of the 1D-BPP in which the bin capacity is 1 and the item sizes $w_i$ are values in $[0,1]$.}

{A function $f: [ 0,1 ] \rightarrow [ 0,1 ]$ is a \emph{Dual Feasible Function} (DFF), if for any {finite set $S$ of non-negative real numbers,} it holds that:
\begin{equation}
\sum_{x \in S} x \leq 1 \implies \sum_{x \in S} f(x) \leq 1.
\end{equation}
Three classes of DFFs were proposed by Fekete and Schepers \cite{FS01}. The survey by
Clautiaux, Alves, and Val{\'e}rio de Carvalho \cite{CAV10} describes relevant (superadditive and maximal) DFFs that can be used in packing contexts. Rietz, Alves, and Val{\'e}rio de Carvalho \cite{RAC10, RAC11, RAC12}, analyzed a number of properties of DFFs (maximality, extremality, worst-case performance) and studied the best parameter tuning for a specific function.}

{Dual feasible functions} can be extended to consider discrete domains and used to derive relaxations for higher-dimensional cutting and packing problems (see Fekete and Schepers \cite{FS04B}).
For the 2-dimensional case, let $f$ and $g$ be two DFFs. Applying $f$ to the item widths and $g$ to the item heights of an instance of any two-dimensional packing problem, we obtain a new instance whose relaxations are valid for the original instance.

Carlier, Clautiaux, and Moukrim \cite{CCM07} proposed DFF-based lower bounds for the 2D-BPP and introduced DFFs depending on the sizes of items and bins ({\em data-dependent} DFFs).
Caprara and Monaci \cite{CM09} solved a bilinear programming problem to determine DFFs that provide the best bound for a given two-dimensional instance.

{Fekete and Schepers \cite{FS04B} extended the concept of DFF to that of conservative scale.
A {\em conservative scale}
is a set of modified item sizes $\widetilde{w}$ (resp. $\widetilde{h}$) such that any subset $\itemset'$ of items {satisfying $\sum_{i \in \itemset'} w_i \leq W$  (resp. $\sum_{i \in \itemset'} h_i \leq H$) also}
satisfies  $\sum_{i \in \itemset'} \widetilde{w}_i \leq W$ (resp. $\sum_{i \in \itemset'} \widetilde{h}_i \leq H$).
Note that a conservative scale is instance dependent while dual feasible functions have general validity. The relationship between (data-dependent) dual feasible functions and conservative scales has been deeply investigated by Belov et al. \cite{BKRS13}.}

{Recently, Serairi and Haouari \cite{SH18} presented a theoretical and experimental study of the most important polynomial-time lower bounding procedures for the 2D-BPP. Among the algorithms they tested, the one by Carlier, Clautiaux, and Moukrim \cite{CCM07} turned out to provide the tightest values.}

{For a comprehensive study on the use of dual feasible functions for integer programming and combinatorial optimization problems the reader is referred to the recent book by Alves et al. \cite{ACCR16}.} \rev{Due to their computational efficiency, dual-feasible functions are often implemented in branch-and-bound algorithms (see, e.g., Clautiaux, Carlier, and Moukrin \cite{CCM07}, and Alvarez-Valdes, Parre{\~n}o and Tamarit \cite{APT09}).}

\subsection{Contiguous One-Dimensional Relaxations} \label{sec:contiguous}
In order to obtain a relaxation of the 2D-SPP, Martello, Monaci, and Vigo \cite{MMV03} introduced the \emph{One-dimensional Contiguous Bin Packing Problem} (1D-CBPP).
The idea is to horizontally ``cut'' each item $i \in \itemset$ into $h_i$ unit-height {\em slices} of width $w_i$ and to solve a 1D-BPP problem with {capacity $W$: the resulting number of one-dimensional bins is then a lower bound on the height of any 2D-SPP solution. To tighten the bound, the 1D-CBPP additionally imposes {\em contiguity}: all slices obtained from an item $i$ must} be packed into $h_i$ {\em consecutive} one-dimensional bins {(see Figure \ref{fig:contiguous_relaxation}(a))}.
\rev{Although the 1D-CBPP is a relaxation of the 2D-SPP, it remains strongly ${\cal NP}$-hard and hence it can only be optimally solved through enumeration. A branch-and-bound algorithm} was presented in \cite{MMV03} for its exact solution and the resulting lower bound was used in an enumerative algorithm for the 2D-SPP.
Later, Alvarez-Valdes, Parre{\~n}o, and Tamarit \cite{APT09} proposed an ILP model for the 1D-CBPP, and used a commercial software to solve it within a branch-and-bound algorithm for the 2D-SPP.

\begin{figure} \centering
\begin{subfigure}[htpb]{0.35\textwidth}\centering
	\begin{tikzpicture}
	\fill[fill=lightgray] (0,0) rectangle (5,6);
	\draw[very thick] (0,0) -- (5,0);
	\draw[very thick] (0,0) -- (0,6);
	\draw[very thick] (5,0) -- (5,6);
	\draw[dash pattern=on 4.78pt off 5pt, line width = 1pt] (0,6) -- (5,6);

	\draw[fill=white,thick] (0,0) rectangle node{1} (2,0.5);
	\draw[fill=white,thick] (0,0.5) rectangle node{1} (2,1);
	\draw[fill=white,thick] (0,1) rectangle node{1} (2,1.5);
	\draw[fill=white,thick] (0,1.5) rectangle node{1} (2,2);
	\draw[fill=white,thick] (0,2) rectangle node{1} (2,2.5);
	\draw[fill=white,thick] (0,2.5) rectangle node{1} (2,3);

	\draw[fill=white,thick] (2,0) rectangle node{2} (3.5,0.5);
	\draw[fill=white,thick] (2,0.5) rectangle node{2} (3.5,1);

	\draw[fill=white,thick] (2,1) rectangle node{3} (4,1.5);
	\draw[fill=white,thick] (2,1.5) rectangle node{3} (4,2);
	\draw[fill=white,thick] (2,2) rectangle node{3} (4,2.5);
	\draw[fill=white,thick] (2,2.5) rectangle node{3} (4,3);
	\draw[fill=white,thick] (1.5,3) rectangle node{3} (3.5,3.5);
	\draw[fill=white,thick] (1.5,3.5) rectangle node{3} (3.5,4);
	\draw[fill=white,thick] (1.5,4) rectangle node{3} (3.5,4.5);

	\draw[fill=white,thick] (3.5,0) rectangle node{4} (4.5,0.5);
	\draw[fill=white,thick] (3.5,0.5) rectangle node{4} (4.5,1);
	\draw[fill=white,thick] (4,1) rectangle node{4} (5,1.5);
	\draw[fill=white,thick] (4,1.5) rectangle node{4} (5,2);

	\draw[fill=white,thick] (4,2) rectangle node{5} (5,2.5);
	\draw[fill=white,thick] (4,2.5) rectangle node{5} (5,3);
	\draw[fill=white,thick] (3.5,3) rectangle node{5} (4.5,3.5);
	\draw[fill=white,thick] (3.5,3.5) rectangle node{5} (4.5,4);
	\draw[fill=white,thick] (3.5,4) rectangle node{5} (4.5,4.5);
	\draw[fill=white,thick] (3,4.5) rectangle node{5} (4,5);

	\draw[fill=white,thick] (0,3) rectangle node{6} (1.5,3.5);
	\draw[fill=white,thick] (0,3.5) rectangle node{6} (1.5,4);
	\draw[fill=white,thick] (0,4) rectangle node{6} (1.5,4.5);
	\draw[fill=white,thick] (0,4.5) rectangle node{6} (1.5,5);

	\draw[fill=white,thick] (1.5,4.5) rectangle node{7} (3,5);

	\draw[dash pattern=on \pgflinewidth off 1.5pt, line width = 1.2pt] (0,0.5) -- (5,0.5);
	\draw[dash pattern=on \pgflinewidth off 1.5pt, line width = 1.2pt] (0,1) -- (5,1);
	\draw[dash pattern=on \pgflinewidth off 1.5pt, line width = 1.2pt] (0,1.5) -- (5,1.5);
	\draw[dash pattern=on \pgflinewidth off 1.5pt, line width = 1.2pt] (0,2) -- (5,2);
	\draw[dash pattern=on \pgflinewidth off 1.5pt, line width = 1.2pt] (0,2.5) -- (5,2.5);
	\draw[dash pattern=on \pgflinewidth off 1.5pt, line width = 1.2pt] (0,3) -- (5,3);
	\draw[dash pattern=on \pgflinewidth off 1.5pt, line width = 1.2pt] (0,3.5) -- (5,3.5);
	\draw[dash pattern=on \pgflinewidth off 1.5pt, line width = 1.2pt] (0,4) -- (5,4);
	\draw[dash pattern=on \pgflinewidth off 1.5pt, line width = 1.2pt] (0,4.5) -- (5,4.5);
	\draw[dash pattern=on \pgflinewidth off 1.5pt, line width = 1.2pt] (0,5) -- (5,5);
	\draw[dash pattern=on \pgflinewidth off 1.5pt, line width = 1.2pt] (0,5.5) -- (5,5.5);

	\end{tikzpicture}
\caption{\label{fig:staircase1}}
\end{subfigure}
\begin{subfigure}[htpb]{0.35\textwidth}\centering
	\begin{tikzpicture}
	\fill[fill=lightgray] (0,0) rectangle (5,6);
	\draw[very thick] (0,0) -- (5,0);
	\draw[very thick] (0,0) -- (0,6);
	\draw[very thick] (5,0) -- (5,6);
	\draw[dash pattern=on 4.78pt off 5pt, line width = 1pt] (0,6) -- (5,6);

	\draw[fill=white,thick] (0,0) rectangle node{1} (2,3);
	\draw[fill=white,thick] (2,0) rectangle node{2} (3.5,1);
	\draw[fill=white,thick] (2,1) rectangle node{3} (4,4.5);
	\draw[fill=white,thick] (4,0) rectangle node{4} (5,2);
	\draw[fill=white,thick] (4,2) rectangle node{5} (5,5);
	\draw[fill=white,thick] (0,3) rectangle node{6} (1.5,5);
	\draw[fill=white,thick] (1.5,4.5) rectangle node{7} (3,5);

	\end{tikzpicture}
\caption{\label{fig:staircase2}}
\end{subfigure}

\caption{(a) {contiguous relaxation of items 1-7}; (b) two-dimensional packing obtained from (a).}
\label{fig:contiguous_relaxation}
\end{figure}
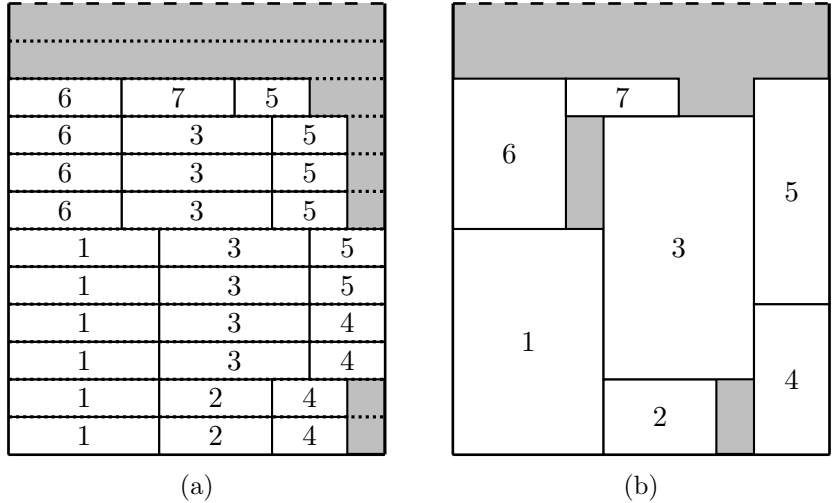

Belov et al. \cite{BKRS09} used the 1D-CBPP as a relaxation of the 2D-OPP. They proposed an alternative ILP model for the 1D-CBPP and a weaker relaxation in which, instead of requiring contiguity,
it is imposed that each one-dimensional bin can contain at most one slice from each item ({\em bar relaxation}). {Friedow and Scheithauer \cite{FS17} solved the 1D-CBPP with a cutting plane approach \rev{based on column generation for} the bar relaxation.} A branch-and-bound algorithm for the 1D-CBPP was later proposed by Mesyagutov et al. \cite{MMBS11}.

C{\^o}t{\'e}, Dell'Amico, and Iori \cite{CDI14} proposed an algorithm for the 2D-SPP based on Benders' decomposition (see Section \ref{sec:exponential} below) of an ILP model.
The resulting master {problem, obtained by ``transposing'' the instance and vertically cutting the
items into unit-width bars, is} a parallel processor scheduling problem with contiguity constraints. The slave problem, \rev{which is still strongly ${\cal NP}$-hard, determines} whether the solution provided by the master can be transformed into a feasible packing {(see Figure \ref{fig:contiguous_relaxation}(b))}. Delorme, Iori, and Martello \cite{DIM17} applied a similar decomposition to solve the 2D-SPP with orthogonal rotations. They proposed an ILP model for the 1D-CBPP based on the classical arc-flow formulation of the 1D-CSP (see Val{\'e}rio de Carvalho \cite{V99}).

\subsection{\rev{State Space Relaxations}}

\rev{The state space relaxation, originally developed by Christofides, Mingozzi, and Toth \cite{CMT81} for routing problems, modifies the state space of a dynamic programming recursion so that its optimal solution is a valid bound for the original problem. This relaxation has been used for the (strongly {\cal NP}-hard) guillotine 2D-KP, for which a state space relaxation can disregard the constraint on the maximum number of copies of each item. The relaxed problem (known in the literature as the {\em unconstrained two-dimensional guillotine knapsack problem}) can be solved in pseudo-polynomial time through dynamic programming (see, e.g., Gilmore and Gomory \cite{GG65}, Beasley \cite{B85g}, Cintra et al. \cite{CMWX08} and Russo et al. \cite{RSS14}).}

\rev{Christofides and Hadjiconstantinou \cite{CH95} strengthened the state space relaxation by associating a non-negative integer multiplier with each item and adopted a subgradient-like procedure to determine good multipliers. Improvements were proposed by Morabito and Pureza \cite{MP10} and by Velasco and Uchoa \cite{VU19}.}

\section{Heuristics}\label{sec:heuristics}
There is a huge literature on approximation algorithms and heuristics for two-dimensional packing problems, and a thorough review of these methods is outside the scope of this survey.
To get an idea of the vastity of this area, it is enough to consider that Ortmann and van Vuuren \cite{OV10} computationally compared a total of 252 heuristics and variants for two-dimensional strip packing problems.
As heuristics are frequently used within exact algorithms (e.g., to initialize the incumbent solution), we provide in the following some pointers to a number of relevant results.

\subsection{Approximation and On-line Algorithms}
Concerning approximation and on-line algorithms, we refer the reader to the {recent (2017) comprehensive} survey by Christensen et al. \cite{CKPT17}, that
gives theoretical insight on fast solution methods characterized by relevant theoretical properties for multi-dimensional packing problems. We only mention here recent relevant works
\rev{that appeared after} the publication of such survey.

G{\'a}lvez et al. \cite{GGHIKW17} proposed a polynomial-time $1.89$-approximation algorithm for the 2D-KP with orthogonal rotations, breaking the previous 2-approximation barrier.

Yu, Mao, and Xiao \cite{YMX16, YMX17b} studied the {\em on-line} 2D-SPP, presenting new lower and upper bounds with guaranteed worst-case performance. In
\cite{YMX17a} they presented an upper bound for the special case of the 2D-SPP with square items. Further improvements were presented by Han et al. \cite{HIYZ16} for the 2D-SPP and by Balogh et al. \cite{BBDEL18}
for the special case of the 2D-BPP in which items are squares.

{For the 2D-SPP, \rev{Henning et al. \cite{HJRF20}} studied pseudo-polynomial approximation algorithms with respect to the
width of the strip, i.e., algorithms whose time complexity is a polynomial function of the strip width. They proved that
there cannot {exist} a pseudo-polynomial algorithm with a ratio better than 5/4 unless ${\cal P} = {\cal NP}$. Jansen and Rau \cite{JR19} closed the gap between inapproximability result and best known algorithm
by presenting an algorithm with approximation ratio 5/4 + $\varepsilon$.}

\subsection{Constructive Heuristics, Local Search, and Metaheuristics}\label{sec:constructive}
A huge number of constructive heuristics for two-dimensional cutting and packing problems can be found in the literature. The reader is referred to Lodi et al. \cite{LMMV14} for classical results in this area. We only mention here two strategies that have been the basis of many solution approaches. The classical \emph{bottom-left algorithm}, proposed in 1980 by Baker, Coffman, and Rivest \cite{BCR80}, packs one item at a time in the  lowest possible position, left justified. {Another classical result (the {\em best-fit} approach)}, proposed by Burke, Kendall, and Whitwell \cite{BKW04}, selects and packs \rev{an} item that better fits the lowest available area.
(See Chazelle \cite{C83} and Imahori and Yagiura \cite{IY10} for efficient implementations of these methods.)

Local search approaches explore the neighborhood of a given solution. Several methods for two-dimensional packing problems are based on fixing part of the solution, unpacking the other items, and completing the solution
according to some strategy. For example, Alvarez-Valdes, Parre{\~n}o, and Tamarit \cite{APT07} create ``holes'' in the current packing and try to use them when completing the solution.
Other approaches pack items according to an input sequence, and perform changes in such sequence to obtain alternative solutions (see, e.g., Burke, Hyde, and Kendall \cite{BHK11} and Wei et al. \cite{WQCX16}).

Constructive heuristics and local search are the base of many efficient metaheuristic algorithms. \rev{We briefly review in the following some results for the two main metaheuristic approaches.}

\subsubsection*{Single-solution Metaheuristics}

Alvarez-Vald{\'e}s, Paraj{\'o}n, and Tamarit \cite{APT02} proposed a GRASP and a Tabu search for the 2D-KP with guillotine {constraints. Alvarez-Valdes, Parre{\~n}o, and Tamarit \cite{APT07, APT08} presented a Tabu search algorithm for the 2D-KP and a GRASP for the 2D-SPP}.
Wei et al. \cite{WOZL11} proposed an effective Tabu search for the two-dimensional rectangular packing problem (see Section \ref{subsec:problems}), both with and without orthogonal rotation. They also embedded it in a
binary search approach to solve the corresponding versions of the 2D-SPP.

Some of the metaheuristics proposed in the literature also work for the {\em three-dimensional} generalization of the 2D-BPP, in which
items and bins are three-dimensional boxes. For both dimensions, Parre{\~ n}o et al. \cite{PAVOT10} presented a combination of GRASP and variable neighborhood search, while Lodi, Martello, and Vigo \cite{LMV04} implemented a unified C computer code based on Tabu search (available at \url{http://or.dei.unibo.it/research_pages/ORcodes/TSpack.html}).

Tabu search heuristics for routing problems with two- and three-dimensional loading/unloading constraints have been proposed by Gendreau et al. \cite{GILM06,GILM08}

\subsubsection*{Population-Based Metaheuristics}
Concerning the 2D-SPP, Iori, Martello, and Monaci \cite{IMM03} developed a hybrid Tabu search-genetic algorithm. Different genetic approaches were presented by Burke et al. \cite{BHKW10}, Matayoshi \cite{M10}, and Borgulya \cite{B14}.

Genetic algorithms for the 2D-KP were proposed by Hadjiconstantinou and Iori \cite{HI07} and by Kierkosz and Luczak \cite{KL14}.
Gon{\c{c}}alves and Resende \cite{GR13} presented a biased random key genetic algorithm for the 2D-BPP and its three-dimensional generalization.

\subsection{Set covering based heuristics}
Monaci and Toth \cite{MT06} developed algorithms for the 2D-BPP based on the set-covering formulation (see Section \ref{sec:column}) induced by a heuristic generation of the patterns.

Cintra et al. \cite{CMWX08} proposed dynamic programming algorithms for the guillotine 2D-KP, both in the $k$-staged and in the non-staged versions. They then used these algorithms in a column generation
heuristic for the corresponding versions of the 2D-CSP and the 2D-SPP. Other heuristics based on column generation for the 2-staged 2D-CSP were presented by Furini et al. \cite{FMDPT12}, by Cui and Zhao \cite{CZ13} {and by Cui, Zhou and Cui \cite{CZC17}}.

\rev{For the 3-staged 2D-CSP, a set covering based heuristic was presented by Vanderbeck \cite{V01}.}

\section{Exact Methods based on Integer Linear Programming Models} \label{sec:ILP}
In this section we review solution approaches that are based on {ILP and {\em Mixed Integer Linear Programming} (MILP)} formulations and hence require the use of a solver. The models may be classified on the basis of their size: polynomial, pseudo-polynomial, or exponential.

\subsection{Polynomial Models}
A mathematical mixed-integer model for a general three-dimensional bin packing problem, involving a polynomial number of variables and {constraints}, was presented by Chen, Lee, and Shen \cite{CLS95}. \rev{(This  obviously implies a polynomial model for the two-dimensional classical problems.)
The model in \cite{CLS95},} that is based on the enumeration of all possible relative placements of each pair of items, can be seen as an extension to the three-dimensional case of the modeling technique proposed by Onodera, Taniguchi, and Tamaru \cite{OTT91} for a two-dimensional block placement problem.

Polynomial formulations for 2-staged guillotine packing were proposed by Lodi and Monaci \cite{LM03} for the 2D-KP and by Lodi, Martello, and Vigo \cite{LMV04b} for the 2D-BPP and the 2D-SPP.
These models were later extended to the 3-staged case of the 2D-BPP by Puchinger and Raidl \cite{PR07} and to the 2-staged 2D-CSP with variable-sized bins by Furini and Malaguti \cite{FM13}.

\subsection{Pseudo-polynomial Models}\label{sec:pseudo}

While polynomial models associate variables to the items, pseudo-polynomial models include variables associated with the positions into the bins where the items {can be packed.}

Beasley \cite{B85} proposed an ILP model for the 2D-KP, based on a discretization of the packing area through a pseudo-polynomial number of two-dimensional {coordinates (see Section \ref{sec:points}) where the bottom-left corner of an item can be packed}. A pseudo-polynomial number of constraints imposes that no unit square is covered by more than one item. This formulation was later used for the solution of other two-dimensional
packing problems. Alvarez-Valdes, Parre{\~n}o, and Tamarit \cite{APT05} adapted the formulation to the special case of the two-dimensional rectangular packing problem (see Section \ref{subsec:problems}) with orthogonal rotations in which all items are identical ({\em Pallet Loading Problem}). {de Queiroz} and Miyazawa \cite{QM14} adapted the model to a variant of the 2D-SPP {in which} the items have to be arranged to form a physically stable packing satisfying a predefined item unloading order. Martello and Monaci \cite{MM15} used a similar model to solve the problems of orthogonally packing a set of rectangles  (with or without allowing rotations) into the smallest square.

Val{\'e}rio de Carvalho \cite{V99} presented an {\em arc-flow} model for the solution of the 1D-CSP. The model is based on a digraph with a pseudo-polynomial number of vertices, in which paths correspond to
feasible packings of a bin and a flow provides an overall problem solution. This formulation was extended to the 2-staged guillotine 2D-CSP by Macedo, Alves, and Val{\'e}rio de Carvalho \cite{MAC10} and to the
2-staged guillotine 2D-SPP by Mrad \cite{M15}.

Another pseudo-polynomial formulation of the 1D-CSP that was further extended to two-dimensional problems is the {{\em one-cut} model, independently} obtained by Rao \cite{R76} and Dyckhoff \cite{D81}.
The model is based on variables that describe the feasible cuts that split a bin (or a residual of a bin) into an item and a residual (or another item).
Silva, Alvelos, and Val{\'e}rio de Carvalho \cite{SAC10} extended this model to the 2-staged and the 3-staged guillotine 2D-CSP, while Furini, Malaguti, and Thomopulos \cite{FMT16} used it for the non-staged \rev{guillotine} 2D-SPP, 2D-CSP, and 2D-KP.

Another family of pseudo-polynomial models derives from time representations of scheduling problems where disjunctive constraints determine the relative positions between each pair of items. Following this methodology, Castro and Oliveira \cite{CO11} developed ILP models for the 2D-SPP, and showed how to extend them to the 2D-BPP and {to the} 2D-KP. Castro and Grossmann \cite{CG12} proposed further improvements for the 2D-SPP.

\subsection{Exponential Models}\label{sec:exponential}
Most exponential models for two-dimensional packing problems are based either on \rev{formulations that associate variables with feasible patterns (see Section \ref{sec:column}), which have an exponential number of variables (columns), or on Benders' decomposition (see below), which produces an exponential number of constraints (rows). The models of the former type can be either set covering/set partitioning models (for bin packing or cutting stock problems) or set packing models (for knapsack problems).}

\subsubsection*{Branch(-and-cut)-and-price}
Models based on the set covering formulation are usually implemented through column generation {(see Section \ref{sec:column})} for solving the associated linear program relaxation. Embedding the linear program into an enumerative scheme allows one to obtain an optimal integer solution through branch-and-price or branch-and-cut-and-price (when cuts are added at the decision nodes).

To the best of our knowledge, Belov and Scheithauer \cite{BS06} were the first to propose a branch-and-cut-and-price algorithm (\rev{based on a set-packing formulation}) for a two-dimensional problem, namely the 2-staged guillotine 2D-KP.

Puchinger and Raidl \cite{PR07} solved the 3-staged guillotine 2D-BPP with a branch-and-price algorithm with special techniques ({\em column generation stabilization}) to accelerate the convergence of the method, as previously suggested by Val{\'e}rio de Carvalho \cite{C05} and Ben Amor, Desrosiers, and Val{\'e}rio de Carvalho \cite{ADC06}.

Pisinger and Sigurd \cite{PS05, PS07} proposed branch-and-price algorithms for {2D-BPPs with variable-sized and fixed-sized} bins. In these algorithms (which are, from a computational point of view, the current state-of-the-art), the columns are generated through decomposition into a one-dimensional knapsack problem and a 2D-OPP, which are {then} solved by constraint programming.
Bettinelli, Ceselli, and Righini \cite{BCR08} developed a branch-and-price approach to the level (2-staged) 2D-SPP in which the pricing problem consists of a penalized one-dimensional knapsack problem. Furini and Malaguti \cite{FM13} \rev{and Mrad, Meftahi and Haouari \cite{MMH13} solved the 2-staged guillotine 2D-CSP, respectively  with
and without variable-sized bins, through branch-and-price algorithms.}

\subsubsection*{Benders' Decomposition} \label{sec:benders}
Already in 1962, Benders \cite{B62} proposed a method to decompose {an MILP} model having a special {\em block structure} into an MILP {\em master problem} and an LP {\em slave problem}. At each iteration, the slave receives the
current master solution and either proves its optimality or generates a cut ({\em Benders' cut}) that is added to the master. The approach is also referred to as {\em row generation}.

Geoffrion \cite{G72} generalized the Benders' decomposition to the case in which the subproblem too is an MILP. The latter decomposition was used by Caprara and Monaci \cite{CM04} and by Baldacci and Boschetti \cite{BB07} to solve the 2D-KP: in both approaches the master is an ILP based on the one-dimensional knapsack problem while the slave solves an associated 2D-OPP.

Hooker and Ottosson \cite{HO03} defined the logic-based Benders' decomposition for general MILP problems: in this approach both the master and the slave are MILPs, but the latter is solved through logical deduction methods, such as constraint programming. This decomposition was successfully applied by Pisinger and Sigurd \cite{PS05, PS07} to the 2D-KP and by {de Queiroz} et al. \cite{QHSM17} to a special 2D-KP that includes pairs of items that cannot both be in the solution (2D-KP with conflict graphs).

Although Benders' decomposition has been successfully applied to cutting and packing problems, the resulting cuts in the master problem may be weak. Codato and Fischetti \cite{CF06} defined stronger  \emph{combinatorial Benders' cuts} for general integer problems. C{\^o}t{\'e}, Dell'Amico, and Iori \cite{CDI14} and Delorme, Iori, and Martello \cite{DIM17} adopted such cuts for the solution of the 2D-SPP without and with item rotation, respectively.
In their logic-based Benders' decomposition, the master problem is a {parallel processor scheduling problem} with contiguity constraints (see Section \ref{sec:contiguous}) that produces the $x$-coordinate of each item, while the slave checks whether feasible $y$-coordinates exist. In the negative case, heuristic methods look for a minimal infeasible set of items that produces a combinatorial Benders' cut.
A similar decomposition, with standard Benders' cuts, was used by C{\^o}t{\'e}, Gendreau, and Potvin \cite{CGP14} for the 2D-OPP with unloading constraints.

\rev{Recently, Martin et al. \cite{MBLMM20} proposed a Benders' decomposition algorithm to solve the guillotine 2D-KP. In their decomposition, the master problem is modeled using the pseudo-polynomial model by Beasley [13], whereas the slave checks the guillotine restriction.}

\section{Exact Methods based on Implicit Enumeration}\label{sec:enumeration}
The algorithms reviewed in this section consist of enumeration schemes that do not explicitly make use of MILP models {and hence do not require the use of a solver}.

\subsection{Branch-and-Bound}\label{sec:B-and-B}
Several enumeration algorithms were derived from the bottom-left strategy proposed by Baker, Coffman, and Rivest \cite{BCR80} {(see Section \ref{sec:constructive}), which produces an} approximate solution by placing one item at a time, in the lowest possible feasible position, left justified. Hadjiconstantinou and Christofides \cite{HC95} developed a tree-search exact algorithm for the 2D-KP
that packs the next item in every possible position such that the item's left and bottom edges touch either the bin or edges of other items ({\em left-most downward} strategy).
Martello and Vigo \cite{MV98} {adopted this strategy in a branch-and-bound algorithm for the 2D-OPP, and embedded it in} a two-level enumeration algorithm for the 2D-BPP.

This strategy later evolved into two classes of implicit enumeration schemes, namely the {\em staircase placement} and the {\em niche placement}.
Both strategies enumerate the possible packing of the items with their bottom-left corner in particular positions ({\em corner points}), induced by the current (partial) packing, whose definition depends on the specific strategy.

\begin{figure} \centering
\begin{subfigure}[htpb]{0.35\textwidth}\centering
	\begin{tikzpicture}
	\draw[fill=lightgray,very thick] (0,0) rectangle (5,5);

	\draw[fill=white,thick] (0,0) rectangle node{1} (1.5,2);
	\draw[fill=white,thick] (1.5,0) rectangle node{2} (3,1.5);
	\draw[fill=white,thick] (3,0) rectangle node{3} (4.5,1);
	\draw[fill=white,thick] (0,2) rectangle node{4} (0.5,4);
	\draw[fill=white,thick] (1.5,1.5) rectangle node{5} (3.5,3);

	\draw[pattern=north west lines] (0.5,2) rectangle (1.5,3);
	\draw[pattern=north west lines] (3,1) rectangle (3.5,1.5);

	\draw[dash pattern=on \pgflinewidth off 1.5pt, line width = 1.4pt] (0,4) -- (0.5,4);
	\draw[dash pattern=on \pgflinewidth off 1.5pt, line width = 1.4pt] (0.5,3) -- (0.5,4);
	\draw[dash pattern=on \pgflinewidth off 1.5pt, line width = 1.4pt] (0.5,3) -- (3.5,3);
	\draw[dash pattern=on \pgflinewidth off 1.5pt, line width = 1.4pt] (3.5,1) -- (3.5,3);
	\draw[dash pattern=on \pgflinewidth off 1.5pt, line width = 1.4pt] (3.5,1) -- (4.5,1);
	\draw[dash pattern=on \pgflinewidth off 1.5pt, line width = 1.4pt] (4.5,0) -- (4.5,1);

	\draw (0,4) node {\small{\textbullet}};
	\draw (0.5,3) node {\small{\textbullet}};
	\draw (3.5,1) node {\small{\textbullet}};
	\draw (4.5,0) node {\small{\textbullet}};
	\end{tikzpicture}
\caption{\label{fig:staircase3}}
\end{subfigure}
\begin{subfigure}[htpb]{0.35\textwidth}\centering
	\begin{tikzpicture}
	\draw[fill=lightgray,very thick] (0,0) rectangle (5,5);

	\draw[fill=white,thick] (0,0) rectangle node{1} (1.5,2);
	\draw[fill=white,thick] (1.5,0) rectangle node{2} (3,1.5);
	\draw[fill=white,thick] (3,0) rectangle node{3} (4.5,1);
	\draw[fill=white,thick] (4.5,0) rectangle node{4} (5,2);
	\draw[fill=white,thick] (1.5,1.5) rectangle node{5} (3.5,3);

	\draw[pattern=north west lines] (3,1) rectangle (4.5,1.5);

	\draw[dash pattern=on \pgflinewidth off 1.5pt, line width = 1.4pt] (0,2) -- (1.5,2);
	\draw[dash pattern=on \pgflinewidth off 1.5pt, line width = 1.4pt] (1.5,2) -- (1.5,3);
	\draw[dash pattern=on \pgflinewidth off 1.5pt, line width = 1.4pt] (1.5,3) -- (3.5,3);
	\draw[dash pattern=on \pgflinewidth off 1.5pt, line width = 1.4pt] (3.5,1.5) -- (3.5,3);
	\draw[dash pattern=on \pgflinewidth off 1.5pt, line width = 1.4pt] (3.5,1.5) -- (4.5,1.5);
	\draw[dash pattern=on \pgflinewidth off 1.5pt, line width = 1.4pt] (4.5,1.5) -- (4.5,2);
	\draw[dash pattern=on \pgflinewidth off 1.5pt, line width = 1.4pt] (4.5,2) -- (5,2);

	\draw (3.5,1.5) node {\small{\textbullet}};

	\draw (0,-0.055) node {};
	\end{tikzpicture}
\caption{\label{fig:staircase4}}
\end{subfigure}
\caption{{corner points of} (a) staircase placement; (b) niche placement.}
\label{fig:placement_strategies}
\end{figure}
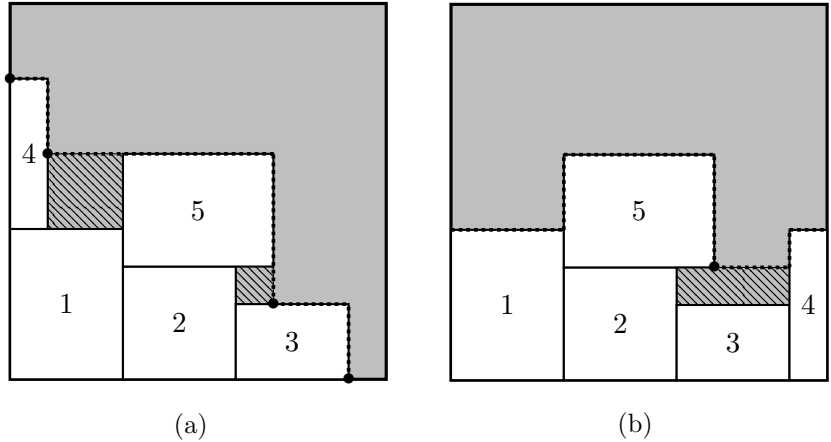

\subsubsection*{Staircase Placement}
The {\em staircase placement} makes use of the {\em envelope} (see Figure \ref{fig:placement_strategies} (a)), defined as the monotone (right-down) staircase-like boundary, that: (i) separates the area where previously packed
items are placed from the area available for the remaining items, and (ii) is composed by segments touching at least one edge of a packed item or the bin. The {\em corner points} {are those in which} the
slope of the envelope changes from vertical to horizontal. These concepts were originally introduced by Scheithauer \cite{SCH97}.

Martello, Monaci, and Vigo \cite{MMV03} proposed a branch-and-bound algorithm for the 2D-SPP which is based on the staircase placement {and on a strategy that avoids the enumeration of} symmetric solutions.
The method was later modified by Iori, Salazar Gonz\'alez, and Vigo \cite{ISGV07} to solve the 2D-OPP with loading/unloading constraints.
Other branch-and-bound algorithms based on {the} staircase placement were proposed by Alvarez-Valdes, Parre{\~n}o, and Tamarit \cite{APT09} for the 2D-SPP,
by Clautiaux, Carlier, and Moukrim \cite{CCM07b} for the 2D-OPP, {and by} Bekrar et al. \cite{BKCS10} for the 2D-SPP with guillotine constraint. The latter algorithm makes use
of a method by Messaoud, Chu, and Espinouse \cite{BCE08} to detect non-guillotine patterns in polynomial time.

{Recently, Xu and Lee \cite{XL18} studied the {{\em continuous berth allocation problem}}, in which incoming vessels need to be assigned a time and a berth location on a quay. Each vessel can be interpreted as a rectangle with given width (the space it occupies in the berth) and height (the time it consumes): the problem is to ``pack'' the vessels into a strip of given width (the total space available at the berth) so as to minimize the total weighted completion times of the activities on the vessels. They solved the problem with a branch-and-bound algorithm that uses a staircase branching rule and obtains at each node a valid lower bound by invoking a model reformulation.}

\subsubsection*{Niche Placement}
The {\em niche placement} uses the {\em skyline structure} (see Figure \ref{fig:placement_strategies} (b)), a boundary composed by orthogonal line segments separating the area of previously packed items from the area that is available for the remaining items. {At each iteration, the {\em niche} is the lowest and (in case of ties) leftmost horizontal segment: the left extreme of the niche is the only {\em corner point} considered for packing the next item.}

\rev{A complete branching is obtained by creating at most $n+1$ nodes as follows. The first $n$ (at most) nodes are created by packing an item in the niche (obviously disregarding
those that would lead to overlapping).  Then, a last node in which no item is packed in the niche is also created. In such node, the niche is closed, so that at the next iteration a new niche will be considered.}

The niche placement was used by Boschetti and Montaletti \cite{BM10} in a branch-and-bound algorithm for the 2D-SPP and by Lesh et al. \cite{LMMM04} for the 2D-SPP with {\em perfect packing} (a packing in which no wasted space is allowed). {Kenmochi et al. \cite{KINYN09} proposed  two branch-and-bound algorithms, respectively adopting niche and staircase placement, for the 2D-SPP with perfect packing, with or without orthogonal rotation.
They additionally derived algorithms for the general 2D-SPP. Improvements of the latter algorithms were presented by Arahori, Imamichi, and Nagamochi \cite{AIN12}.}

\subsubsection*{Other Enumeration Schemes}

{Certain solution methods for the 2D-KP are based on an enumeration of candidate subsets of items {combined with} an inner solution method that checks whether the current subset can be feasibly packed.
This method has been used by Caprara and Monaci \cite{CM04} for the 2D-KP and by Fekete, Schepers, and van der Veen \cite{FSV07} for knapsack problems with an arbitrary number of dimensions.}

{Christofides and Whitlock \cite{CW77} and Christofides and Hadjiconstantinou \cite{CH95} introduced enumerative approaches for the 2D-KP with guillotine constraints. Their algorithms are based on a tree-search approach
where branchings correspond to cuts on a rectangle and bounds are obtained by solving relaxations through dynamic programming.
Dolatabadi, Lodi, and Monaci \cite{DLM12} improved this enumeration scheme, obtaining a branch-and-bound and a branch-and-cut algorithm.} \rev{Another search strategy that has been extensively used for problems with guillotine constraints is the bottom-up approach (see, e.g., Viswanathan and Bagchi \cite{VB93}, Hifi \cite{H97}, Cung, Hifi, and Le Cun \cite{CHLC00}, and Fleszar \cite{F16}). {It is based on enumerating all feasible patterns obtained through {\em builds}: a build is a rectangle produced by the (horizontal or vertical) combination of a pair of items or other builds.} The bottom-up approach also inspired the development of MILP formulations for guillotine problems (see, e.g., Martin, Morabito, and Munari \cite{MMM20}).}

{A special case of 2-staged guillotine cutting patterns are the {\em checkerboard patterns}. While the guillotine constraint allows each cut to separate a bin into two pieces, {which} can then be treated as new (smaller) bins, the checkerboard constraint imposes that the solution be obtained through a set of horizontal cuts and a set of vertical cuts, all performed on the original bin. Yanasse and Katsurayama \cite{YK05, YK08} presented enumerative algorithms for the 2D-KP with checkerboard patterns.
}

\rev{Clautiaux et al. \cite{CSVV18} proposed label setting algorithms for the four-staged guillotine 2D-KP. They presented reduction procedures, filtering rules based on a Lagrangian relaxation, and a state space relaxation.}

\subsection{{Graph-Based Approaches}}

\rev{Morabito and Arenales \cite{MA96} proposed an exact algorithm for the guillotine 2D-KP based on AND/OR graphs. In an AND/OR-graph representation, the nodes correspond to rectangles, the arcs correspond to cuts, and
cutting patterns are represented as {\em complete paths} in the graph.}

{Fekete and Schepers \cite{FS04A} proposed a graph-theoretical characterization of the packing of a set of items into a bin. Their representation makes use of two interval graphs, associated with the horizontal and vertical {dimensions}, respectively: each item corresponds to a vertex, and two vertices are connected by an edge if and only if the projections of the corresponding items on the horizontal/vertical axis
overlap. Theoretical properties allow one to detect the feasibility of a packing. This {\em interval graph model} easily extends to packings in higher dimensions.}

{Based on the interval graph representation, Fekete, Schepers, and van der Veen \cite{FSV07} proposed a branch-and-bound algorithm for higher-dimensional orthogonal packing problems. The algorithm was later improved
by Belov and Rohling \cite{BR13} through LP bounds based on the bar relaxation (see Section \ref{sec:contiguous}). In order to avoid the enumeration of symmetrical interval graphs and
unnecessarily enumerated packing classes (as observed by Ferreira and Oliveira \cite{FO08}), Joncour and P{\^e}cher \cite{JP10} proposed to adopt the so-called
{\em consecutive ones matrices} to enumerate relevant interval graphs. Further improvements for the 2-dimensional case were proposed by Joncour, P{\^e}cher, and Valicov \cite{JPV12}.}

Clautiaux, Jougler, and Moukrim \cite{CJM08} presented a graph-theoretical model for the 2D-OPP with guillotine constraints. The model is based on an arc-colored directed graph, which can be replaced by an uncolored undirected multigraph, and can be used for solving the problem through \rev{constraint programming}.

\subsection{Constraint Programming}

{As previously mentioned (see Section \ref{sec:exponential}), Pisinger and Sigurd \cite{PS05, PS07} used a constraint programming formulation to solve 2D-OPP instances arising as subproblems in their decomposition approaches
for the 2D-BPP (with {fixed-sized or variable-sized} bins). The main constraints are {related to the} relative positioning of each pair of items.  This strategy was later adopted for other two-dimensional packing problems, e.g., by Korf, Moffitt, and Pollack \cite{KMP10} for the problem of orthogonally packing a set of two-dimensional items into a rectangle having minimum area, and by {de Queiroz} et al. \cite{QHSM17} for a variant of the 2D-KP (see Section \ref{sec:benders}).}

{Clautiaux et al. \cite{CJCM08} proposed a constraint-based scheduling model for the 2D-OPP, and solved it through constraint programming and effective propagation techniques. The approach was improved by Mesyagutov, Scheithauer, and Belov \cite{MSB12} in their algorithms for the 2D-SPP and the 2D-OPP, by embedding LP-based pruning rules into the constraint propagation process.}

{Soh et al. \cite{SITBN10} solved the 2D-OPP by {iteratively reducing it to {\em satisfiability testing} (SAT) problems, which look for a feasible assignment of} a set of boolean variables. The approach extends to the solution of the 2D-SPP. Grandcolas and Pinto \cite{GP10} proposed a SAT encoding of the interval graph model by Fekete and Schepers \cite{FS04A} for higher dimensional problems, and compared its efficiency with that of other SAT encodings.}

{Delorme, Iori, and Martello \cite{DIM17} {proposed, for the 2D-CSP, a} constraint programming formulation based on non-overlapping intervals to determine whether a solution of the contiguous one-dimensional relaxation (see Section \ref{sec:contiguous}) produces a feasible two-dimensional packing.}

\rev{
\section{Open Problems} \label{sec:open_problems}
}
\rev{
Our study shows that considerable improvements in the exact solution of two-dimensional cutting and packing problems emerged in the last decades, typically allowing to determine an optimal solution even for large-size instances. However, the inherent hardness of this class of problems is witnessed by the existence of unsolved instances with relatively few items. In the following, we provide a list of the main benchmarks that are still unsolved to proven optimality, not only for the main problems addressed in this survey, but also for other relevant problem variants:
\begin{itemize}
\item for the 2D-SPP, specific instances with 20 items cannot be solved to proven optimality even by recent and sophisticated algorithms (see, e.g., C{\^o}t{\'e}, Dell'Amico, and Iori \cite{CDI14}). Among the 500 instances of the so called {\em 10 classes} benchmark (proposed a long time ago by Berkey and Wang \cite{BW87} and Martello and Vigo \cite{MV98} for the 2D-BPP), only 322 are solved to proven optimality. A number of instances with 40 items are still open for the 2D-SPP with guillotine cuts (see, e.g., Mrad [142]);
\item the 2D-SPP with item rotation is even more difficult from a computational perspective: among the above mentioned 500 instances, just 176 could be solved to optimality,
and in particular 56 instances with only 20 items are still open (see, e.g., Delorme, Iori, and Martello \cite{DIM17});
\item the work by Delorme, Iori, and Martello \cite{DIM17} also presents extensive computational results for the {\em pallet loading problem} (a 2D-OPP variant in which all items have the same dimensions and rotation of 90 degrees is allowed) and the {\em rectangle packing problem} (a 2D-SPP variant in which items must be packed, without rotation, in a square of minimum area, see Martello and Monaci \cite{MM15}). For the former problem, a few dozen instances remain unsolved (see \url{http://lagrange.ime.usp.br/~lobato/packing/cover3.php} for an up-to-date list). For the latter problem, more than 200 instances are unsolved, 50 of these (set RND\_R15
    in \cite{MM15}) involving just 15 items;
\item after more than 20 years, many 2D-BPP instances of the 10 classes benchmark (see \cite{BW87}, \cite{MV98}) with $n \geq 60$  are still unsolved, and the same holds even for some instances with $n = 40$ (see Pisinger and Sigurd \cite{PS07}). For all such instances the difference between the best upper and lower bounds is one bin;
\item Pisinger and Sigurd \cite{PS05} created 500 instances of the 2D-BPP with variable bin sizes and costs, by modifying the 10 classes above. They could only solve 164 of them, leaving as open problems 29 instances with 20 items, and 62 with 40 items. All instances are available at \url{https://www.computational-logistics.org/orlib/topic/2D\%20Variable-sized\%20Bin\%20Packing/index.html\#intro~} ;
\item a very difficult 2D-KP instance with only 32 items, proposed 35 years ago by Beasley \cite{B85}, is still open (see Caprara and Monaci \cite{CM04} and Baldacci and Boschetti \cite{BB07}). For the 2D-KP with guillotine constraints, several   instances with either 25 or 50 items, recently proposed by Velasco and Uchoa \cite{VU19}, are open;
\item for the 2D-OPP, two difficult sets of instances have been created by Mesyagutov, Scheithauer, and Belov \cite{MSB12}, and can be downloaded at \url{http://www.math.tu-dresden.de/~capad/TESTS/OPP/12_cp2_data.zip}. The zipfile includes, among others, 1080 instances ({\tt Gleb\_lpcs} and {\tt Gleb\_opp\_gen}) with 20 items and bin width equal to either 100 or 1000, 18 of them being still unsolved. Other sets, known as C, N, and T (available at \url{https://www.euro-online.org/websites/esicup/data-sets/\#1535972088188-55fb7640-4228}), contain 91 instances with about 70 items and bin width at most 200: only 39 of them have been solved to optimality (see C{\^o}t{\'e} and Iori \cite{CI18});
\item for the 2D-OPP with unloading constraints, C{\^o}t{\'e}, Gendreau, and Potvin \cite{CGP14} generated 6 classes containing in total 3282 instances, with 26 items per instance on average: only 2179 of them have been solved to proven optimality. The instances can be downloaded from \url{https://w1.cirrelt.ca/~cotejean/} (see 2OPP-UL).
\end{itemize}
Most of the above instances can also be downloaded from Internet repositories like, e.g., the well-known {\em OR Library} \url{http://people.brunel.ac.uk/~mastjjb/jeb/info.html} by John Beasley \cite{B90} or the library of the OR group of the University of Bologna \url{http://or.dei.unibo.it/library}. {In addition, a generator of instances for two-dimensional rectangular cutting and packing problems, proposed by Silva, Oliveira, and W{\"a}scher \cite{SOW14}, can be \rev{downloaded} from \url{https://sites.google.com/gcloud.fe.up.pt/cutting-and-packing-tools}.}
}

\section{{Conclusions \rev{and Future Research Directions}}}\label{sec:conclusions}
\rev{We reviewed over 180 papers related to two-dimensional orthogonal cutting and packing problems. The literature in this field has grown considerably in recent years,
as also shown by the fact that about half of our references appeared in the last ten {years, and just 20\% were published before 2000}.}
We described the main preprocessing and relaxation methods and briefly examined heuristic and approximation algorithms. We discussed mathematical models and reviewed the most effective implicit enumeration approaches. \rev{Finally, we provided an extensive list of instances for which the optimal solution is still unknown.}

\rev{To facilitate future research, we provide in Table \ref{tab:table1} a summary of the publications from the last 20 years on exact methods for 2D problems. The table presents the main problems, variants and techniques studied by each paper. The problems are divided into the four main categories discussed in Section \ref{sec:problems} (column BPP gives both references to 2D-BPP and 2D-CSP). The variants are given in the columns OR (orthogonal rotation), GL (guillotine cuts), VS (variable sized bins), LU (loading/unloading constraints) and OTH (other variants). The main techniques are shown in columns MILP, B\&B (branch-and-bound), B\&P (branch-(-and-cut)-and-price), BD (Benders' decomposition), GB (graph based), CP (constraint programming), and REL (relaxations).}

\rev{Some statistics on the recent exact methods can be drawn from the references in the table: 20 papers studied the 2D-SPP, 17 studied the 2D-BPP and/or the 2D-CSP, 14 studied the 2D-KP, and 19 studied the 2D-OPP. The most studied variant is the guillotine case. The most popular solution methods are MILP and branch-and-bound (both with 14 references). Note, however, that the number of papers on Benders' decomposition is relatively high although the use of this technique for cutting and packing problems is quite recent.}

\rev{
We also provide the reader some hints of what we consider promising future research directions:
\begin{itemize}
\item Constraint Programming, either as a stand-alone algorithm or as a sub-problem in decomposition methods, recently led to consistent improvements in computational results on 2D problems. Further research on this type of methodology is envisaged, as we see the potential to obtain further improvements;
\item the use of Benders' decompositions (with or without Constraint Programming) for solving 2D problems is quite recent. For some of the classical problems, these decompositions represent the state-of-the-art, as they could solve to proven optimality several open instances. Such techniques are easily adaptable to many other problem variants, as one can consider embedding the additional features either in the master or in the slave. We thus envisage further research on this decomposition to solve new problem variants;
\item current research trends appear to be intensive for 2D problems with additional constraints, like loading/unloading, stability, fragility, etc.;
\item most of the techniques we presented can be generalized to the case of three or more dimensions, although this is not always straightforward. The amount of research on three-dimensional cutting and packing problems has been so far quite limited, in spite of their large number of applications. Interesting research directions can be explored for these problems, starting from the 2D ones.
\end{itemize}
}

{We hope that this review will encourage researchers to pursue investigations in these fascinating topics where there is still room for improving the algorithmic approaches.}

\begin{table} \caption{\rev{Summary of papers with exact methods and relaxations from the last 20 years.}} \label{tab:table1}
\rev{
\footnotesize
\renewcommand{\arraystretch}{0.85}
\setlength{\tabcolsep}{1.8mm}
\rowcolors{2}{white}{gray!15}
\begin{tabular}{@{}l|cccc|ccccc|ccccccc@{}}
\toprule
\multicolumn{1}{c}{}                                                                          & \multicolumn{4}{c}{2D Problems} & \multicolumn{5}{c}{Variants} & \multicolumn{7}{c}{Techniques}          \\*
\cmidrule(lr){2-5} \cmidrule(lr){6-10} \cmidrule(lr){11-17}
\rowcolor{white}
\multicolumn{1}{c}{Reference }                                                                                                & \RotText{SPP}   & \RotText{BPP}   & \RotText{KP}   & \RotText{OPP}   & \RotText{OR}  & \RotText{GL}  & \RotText{VS}  & \RotText{LU}  & \RotText{OTH}  & \RotText{MILP} & \RotText{B\&B} & \RotText{B\&P} & \RotText{BD} & \RotText{GB} & \RotText{CP} & \RotText{REL} \\*
\cmidrule(lr){1-1} \cmidrule(lr){2-5} \cmidrule(lr){6-10} \cmidrule(lr){11-17}
Alvarez-Valdes, Parre{\~n}o, and Tamarit \cite{APT09}       					                   &       &       &      &       &     &     &     &     &      &      & \checkmark    &      &    &    &    & \checkmark   \\
Arahori, Imamichi, and Nagamochi\cite{AIN12}                                                       & \checkmark     &       &      &       &     &     &     &     &      &      & \checkmark    &      &    &    &    & \checkmark   \\
Baldacci and Boschetti \cite{BB07}                                                                 &       &       & \checkmark    &       &     &     &     &     &      &      &      &      & \checkmark  &    &    &     \\
Bekrar et al. \cite{BKCS10}                                                                        & \checkmark     & \checkmark     &      &       &     & \checkmark   &     &     &      &      & \checkmark    &      &    &    &    & \checkmark   \\
Belov et al. \cite{BKRS09}                                                                         &       &       &      & \checkmark     &     &     &     &     &      &      &      &      &    &    &    & \checkmark   \\
Belov et al. \cite{BKRS13}                                                                         &       &       &      & \checkmark     &     &     &     &     &      &      &      &      &    &    &    & \checkmark   \\
Belov and Rohling \cite{BR13}                                                                      &       &       &      & \checkmark     &     &     &     &     &      &      &      &      &    & \checkmark  &    & \checkmark   \\
Belov and Scheithauer \cite{BS06}                                                                  &       &       & \checkmark    &       &     & \checkmark   &     &     &      &      &      & \checkmark    &    &    &    &     \\
Bettinelli, Ceselli, and Righini \cite{BCR08}                                                      & \checkmark     &       &      &       &     & \checkmark   &     &     &      &      &      & \checkmark    &    &    &    &     \\
Boschetti, Hadjiconstantinou, and Mingozzi \cite{BHM02}                                            &       &       & \checkmark    &       &     &     &     &     &      & \checkmark    &      &      &    &    &    & \checkmark   \\
Boschetti and Mingozzi \cite{BM03A}                                                                &       & \checkmark     &      &       &     &     &     &     &      &      &      &      &    &    &    & \checkmark   \\
Boschetti and Mingozzi \cite{BM03B}                                                                &       & \checkmark     &      &       & \checkmark   &     &     &     &      &      &      &      &    &    &    & \checkmark   \\
Boschetti and Montaletti \cite{BM10}                                                              & \checkmark     &       &      &       &     &     &     &     &      &      & \checkmark    &      &    &    &    & \checkmark   \\
Caprara and Monaci \cite{CM04}                                                                     &       &       & \checkmark    &       &     &     &     &     &      &      &      &      & \checkmark  &    &    &     \\
Caprara and Monaci \cite{CM09}                                                                     &       & \checkmark     &      &       &     &     &     &     &      &      &      &      &    &    &    & \checkmark   \\
Carlier, Clautiaux, and Moukrim \cite{CCM07}                                                       &       & \checkmark     &      &       &     &     &     &     &      &      &      &      &    &    &    & \checkmark   \\
Castro and Grossmann \cite{CG12}                                                                   & \checkmark     &       &      &       &     &     &     &     &      & \checkmark    &      &      &    &    &    &     \\
Castro and Oliveira \cite{CO11}                                                                    & \checkmark     & \checkmark     & \checkmark    &       &     &     &     &     &      & \checkmark    &      &      &    &    &    &     \\
Cintra et al. \cite{CMWX08}                                                                        & \checkmark     & \checkmark     & \checkmark    &       & \checkmark   & \checkmark   &     &     &      &      &      &      &    &    &    & \checkmark   \\
Clautiaux, Carlier, and Moukrim \cite{CCM07}                                                       &       &       &      & \checkmark     &     &     &     &     &      &      & \checkmark    &      &    &    &    & \checkmark   \\
Clautiaux et al.\cite{CJCM08}                                                                      &       &       &      & \checkmark     &     &     &     &     &      &      &      &      &    &    & \checkmark  &     \\
Clautiaux, Jouglet, and Moukrim \cite{CJM08}                                                       &       &       &      & \checkmark     &     & \checkmark   &     &     &      &      &      &      &    & \checkmark  &    &     \\
Clautiaux et al. \cite{CSVV18}                                                                     &       &       & \checkmark    &       &     & \checkmark   &     &     &      &      & \checkmark    &      &    &    &    & \checkmark   \\
C{\^o}t{\'e}, Dell'Amico, and Iori \cite{CDI14} 												   & \checkmark     &       &      & \checkmark     &     &     &     &     &      &      &      &      & \checkmark  &    &    &     \\
C{\^o}t{\'e}, Gendreau, and Potvin \cite{CGP14} 												   &       &       &      & \checkmark     &     &     &     & \checkmark   &      &      &      &      & \checkmark  &    &    &     \\
Delorme, Iori, and Martello \cite{DIM17}                                                           & \checkmark     &       &      &       & \checkmark   &     &     &     &      &      &      &      & \checkmark  &    &    &     \\
Dolatabadi, Lodi, and Monaci \cite{DLM12}                                                          &       &       & \checkmark    &       &     & \checkmark   &     &     &      &      & \checkmark    &      & \checkmark  &    &    &     \\
Fekete and Schepers \cite{FS04A}                                                                   &       &       &      & \checkmark     &     &     &     &     &      &      &      &      &    & \checkmark  &    &     \\
Fekete and Schepers \cite{FS04B}                                                                   &       &       &      & \checkmark     &     &     &     &     &      &      &      &      &    &    &    & \checkmark   \\
Fekete, Schepers, and Veen \cite{FSV07}                                                            &       &       &      & \checkmark     &     &     &     &     &      &      & \checkmark    &      &    & \checkmark  &    &     \\
Fleszar \cite{F16}                                                                                 &       &       &      & \checkmark     & \checkmark   & \checkmark   &     &     &      &      & \checkmark    &      &    &    &    &     \\
Friedow and Scheithauer \cite{FS17}                                                                & \checkmark     &       &      &       &     &     &     &     &      &      &      &      &    &    &    & \checkmark   \\
Furini and Malaguti \cite{FM13}                                                                    &       & \checkmark     &      &       &     & \checkmark   & \checkmark   &     &      & \checkmark    &      & \checkmark    &    &    &    &     \\
Furini, Malaguti, and Thomopulos \cite{FMT16}                                                      & \checkmark     & \checkmark     & \checkmark    &       &     & \checkmark   &     &     &      & \checkmark    &      &      &    &    &    &     \\
Grandcolas and Pinto \cite{GP10}                                                                   &       &       &      & \checkmark     &     &     &     &     &      &      &      &      &    &    & \checkmark  &     \\
Iori, Salazar-Gonzalez, and Vigo \cite{ISV07}                                                      &       &       &      & \checkmark     &     &     &     & \checkmark   &      &      & \checkmark    &      &    &    &    &     \\
Joncour and P{\^e}cher \cite{JP10}                                                                     &       &       &      & \checkmark     &     &     &     &     &      &      &      &      &    & \checkmark  &    &     \\
Joncour, P{\^e}cher, and Valicov \cite{JPV12}                                                          &       &       &      & \checkmark     &     &     &     &     &      &      &      &      &    & \checkmark  &    &     \\
Kenmochi et al. \cite{KINYN09}                                                                     & \checkmark     &       &      &       & \checkmark   &     &     &     &      &      & \checkmark    &      &    &    &    &     \\
Lesh et al. \cite{LMMM04}                                                                          & \checkmark     &       &      &       &     &     &     &     &      &      & \checkmark    &      &    &    &    &     \\
Lodi, Martello, and Vigo \cite{LMV04b}                                                             & \checkmark     & \checkmark     &      &       &     & \checkmark   &     &     &      & \checkmark    &      &      &    &    &    & \checkmark   \\
Lodi and Monaci \cite{LM03}                                                                        &       &       & \checkmark    &       &     & \checkmark   &     &     &      & \checkmark    &      &      &    &    &    &     \\
Macedo, Alves, and Val{\'e}rio de Carvalho \cite{MAC10} 				                           &       & \checkmark     &      &       &     & \checkmark   &     &     &      & \checkmark    &      &      &    &    &    &     \\
Martello, Monaci, and Vigo \cite{MMV03}                                                            & \checkmark     &       &      &       &     &     &     &     &      &      & \checkmark    &      &    &    &    & \checkmark   \\
Martin et al. \cite{MBLMM20}                                                                       &       &       & \checkmark    &       &     & \checkmark   &     &     &      &      &      &      & \checkmark  &    &    &     \\
Martin, Morabito, and Munari \cite{MMM20}                                                                       &       &       & \checkmark    &       &     & \checkmark   &     &     &      &  \checkmark &      &      & &    &    &     \\
Matayoshi \cite{M10}                                                                               &       &       &      &       &     &     &     &     &      &      &      &      &    &    &    &     \\
Melega, de Araujo, and Jans \cite{MAJ18}                                                           &       &       &      &       &     &     &     &     &      &      &      &      &    &    &    &     \\
Messaoud et al. \cite{BCE08}                                                                            & \checkmark     & \checkmark     &      & \checkmark     &     & \checkmark   &     &     &      & \checkmark    &      &      &    &    &    &     \\
Mesyagutov, Scheithauer, and Belov \cite{MSB12}                                                    & \checkmark     &       &      & \checkmark     &     &     &     &     &      &      &      &      &    &    & \checkmark  & \checkmark   \\
Mesyagutov et al. \cite{MMBS11}                                                                    &       &       &      &       &     &     &     &     &      &      &      &      &    &    &    & \checkmark   \\
Mrad \cite{M15}                                                                                    & \checkmark     &       &      &       &     & \checkmark   &     &     &      & \checkmark    &      &      &    &    &    &     \\
Pisinger and Sigurd \cite{PS05}                                                                    &       & \checkmark     &      &       &     &     & \checkmark   &     &      &      &      & \checkmark    &    &    &    &     \\
Pisinger and Sigurd \cite{PS07}                                                                    &       & \checkmark     &      &       &     &     &     &     &      &      &      & \checkmark    &    &    &    &     \\
Puchinger and Raidl \cite{PR07}                                                                    &       & \checkmark     &      &       &     & \checkmark   &     &     &      & \checkmark    &      & \checkmark    &    &    &    &     \\
Queiroz et al. \cite{QHSM17}                                                                       &       &       & \checkmark    &       &     &     &     &     & \checkmark    &      &      &      & \checkmark  &    &    &     \\
Queiroz and Miyazawa \cite{QM14}                                                                   & \checkmark     &       &      &       &     &     &     &     & \checkmark    & \checkmark    &      &      & \checkmark  &    &    &     \\
Serairi and  Haouari \cite{SH18}                                                                   &       & \checkmark     &      &       &     &     &     &     &      &      &      &      &    &    &    & \checkmark   \\
Silva, Alvelos, and Val{\'e}rio de Carvalho \cite{SAC10}                				           &       & \checkmark     &      &       &     & \checkmark   &     &     &      & \checkmark    &      &      &    &    &    &     \\
Soh et al. \cite{SITBN10}                                                                           & \checkmark     &       &      &       &     &     &     &     &      &      &      &      &    &    & \checkmark  &     \\
Velasco and Uchoa \cite{VU19}                                                                      &       &       & \checkmark    &       &     & \checkmark   &     &     &      &      &      &      &    &    &    & \checkmark   \\
Yanasse and Katsurayama \cite{YK05}                                                                &       &       &      & \checkmark     &     &     &     &     & \checkmark    &      & \checkmark    &      &    &    &    &     \\* \bottomrule
\end{tabular}
}
\end{table}

\section*{Acknowledgments}
{Research supported by Air Force Office of Scientific Research, by CNPq (under grant numbers 425340/2016-3 and 314366/2018-0), and by FAPESP (under grant numbers 2015/11937-9, 2016/01860-1, 2016/23552-7, 2018/19217-3, and 2019/12728-5). 

\bibliographystyle{plain}
\bibliography{2D-Survey-Rev}

\end{document}